\documentclass{elsarticle}
\usepackage{a4wide}
\usepackage[english]{babel}
\usepackage{amsmath}
\usepackage{amssymb}
\usepackage{siunitx}
\usepackage{amsfonts}
\usepackage{subcaption}
\usepackage{xcolor}
\usepackage[normalem]{ulem}
\usepackage{cancel}
\usepackage{textcomp}

\biboptions{sort&compress}

\journal{Computer Methods in Applied Mechanics and Engineering}

\newtheorem{remark}{Remark}[section]

\graphicspath{{images/}}

\begin{document}

\begin{frontmatter}

\title{Preconditioning immersed isogeometric finite element methods with application to flow problems}

\author[TUE]{F. de Prenter\corref{mycorrespondingauthor}}
\cortext[mycorrespondingauthor]{Corresponding author. Tel.: +31 61 516 2599}
\ead{f.d.prenter@tue.nl}

\author[TUE]{C.V. Verhoosel}
\ead{c.v.verhoosel@tue.nl}

\author[TUE]{E.H. van Brummelen}
\ead{e.h.v.brummelen@tue.nl}

\address[TUE]{Department of Mechanical Engineering, Eindhoven University of Technology, The Netherlands}

\begin{abstract}
Immersed finite element methods generally suffer from conditioning problems when cut elements intersect the physical domain only on a small fraction of their volume. De Prenter \emph{et al.}\ [Computer Methods in Applied Mechanics and Engineering, 316 (2017) pp. 297--327] present an analysis for symmetric positive definite (SPD) immersed problems, and for this class of problems an algebraic preconditioner is developed. In this contribution the conditioning analysis is extended to immersed finite element methods for systems that are not SPD and the preconditioning technique is generalized to a connectivity-based preconditioner inspired by Additive-Schwarz preconditioning. This Connectivity-based Additive-Schwarz (CbAS) preconditioner is applicable to problems that are not SPD and to mixed problems, such as the Stokes and Navier-Stokes equations. A detailed numerical investigation of the effectivity of the CbAS preconditioner to a range of flow problems is presented.
\end{abstract}


\end{frontmatter}


\section{Introduction}
Immersed finite element methods have been demonstrated to have great potential for problems that are posed on domains for which traditional (mesh-fitting) meshing techniques encounter problems, such as prohibitively large meshing costs or the necessity for manual intervention. In recent years, immersed methods -- such as the finite cell method \cite{Parvizian2007}, CutFEM \cite{Burman2015} and immersogeometric analysis \cite{Kamensky2015} -- have been a valuable companion to isogeometric analysis \cite{Hughes2005} as they enable computations on volumetric domains based on the availability of merely a CAD surface representation \cite{Schillinger2012,Rank2012} or voxelized geometries \cite{Schillinger2015}. Additionally, in the context of isogeometric analysis, immersed techniques can be considered as a natural way to incorporate CAD trimming curves in the design-through-analysis cycle \cite{Ruess2013,Ruess2014,Guo2015,Guo2017}.

A disadvantage of immersed finite element methods is that they can result in severely ill-conditioned matrices when the system contains elements that only intersect the physical domain on a small fraction of their volume, see, \emph{e.g.,}\ \cite{Burman2010,Sanches2011,BurmanHansbo2012,Massing2014,Rueberg2014,Schillinger2015,Dettmer2016,Lehrenfeld2016,SIPIC}. We note that ill-conditioning effects due to small volume fractions are not exclusive to finite element methods, and also occur for immersed finite volume methods, see \emph{e.g.,}\ \cite{Johansen1998}. In \cite{SIPIC} we have analyzed these conditioning problems for partial differential equations that lead to symmetric and coercive variational forms. In particular we considered isogeometric discretizations of the Poisson problem and linear elasticity, both with symmetrically imposed Dirichlet conditions by means of Nitsche's method \cite{Nitsche1971}. Based on this analysis the SIPIC (\emph{Symmetric Incomplete Permuted Inverse Cholesky}) preconditioner was developed as an algebraic preconditioning technique tailored to the above-mentioned class of partial differential equations.

The SIPIC preconditioner exhibits excellent behavior for symmetric and coercive variational forms -- and therefore to Symmetric Positive Definite (SPD) matrices -- which is the class of problems to which it is restricted. Besides this restriction to SPD problems, a drawback of SIPIC is that its robustness can deteriorate for high-order bases with low regularity such as $p$-FEM bases (\emph{e.g.,}\ \cite{Duester2017}).

In the present work we show that SIPIC can be interpreted as a sparse Additive-Schwarz preconditioner (\emph{e.g.,}\ \cite{Smith1996,Toselli2005}) and present a generalization exploiting the connectivity data of the basis. This \emph{Connectivity-based Additive-Schwarz} (CbAS) preconditioner does not suffer from the restrictions of SIPIC and is robust for high-order bases with low regularity. Additionally, we present a method based on the Schur complement to efficiently apply the CbAS preconditioner to mixed variational forms. These developments broaden the range of applications and, most notably, open the doors to applications in flow problems. Flow problems on immersed grids have been studied for decades, see \emph{e.g.,}\ the pioneering work in \cite{Peskin1972,Peskin1977} and the more recent review article \cite{Peskin2002}. Recent work on immersed flow problems involves numerous applications, such as: flows around complex (CAD) objects \cite{Bazilevs2012,Xu2015,Hsu2016}; fluid-structure interaction with large deformations \cite{Burman2014,Massing2014,Wang2017,Wu2017,Schott2014,Schott2016,Kadapa2017}; multiphase flows \cite{Schott2015}; topology optimization \cite{Villanueva2017} and flow problems on scanned domains such as, \emph{e.g.,}\ the imbibition of porous media or biomechanical applications \cite{Hsu2014,Kamensky2015,Hsu2015}. However, robust preconditioning of such problems is still an open research topic. We note that alternative broadly applicable approaches to resolving the ill-conditioning of immersed problems exist, of which the ghost penalty \cite{Burman2010,BurmanHansbo2012} -- which is customary for methods referred to as CutFEM \cite{Burman2015} -- and basis function manipulation \cite{Hoellig2001,Hoellig2005,Rueberg2012,Rueberg2014,Rueberg2016} are particular noteworthy. Also alternative preconditioners have been developed, such as \cite{Menk2011,Lehrenfeld2016,Badia2017}. Similar to SIPIC these preconditioners are, however, restricted to linear elasticity and the Poisson problem however.

The CbAS preconditioner developed herein is tested for a range of problems of increasing complexity: a Poisson problem with boundary conditions imposed by the symmetric and nonsymmetric Nitsche method; an SUPG-stabilized convection-dominated convection-diffusion problem with the (standard) symmetric Nitsche method; a Stokes flow with the symmetric Nitsche method; and the steady and transient incompressible Navier-Stokes equations, also with the symmetric Nitsche method. These four problems enable investigation of the conditioning of symmetric and nonsymmetric single-field and two-field immersed problems. The symmetric and nonsymmetric Nitsche methods for the Poisson problem furthermore enable a comparison between the conditioning properties of both methods.

In Section~\ref{sec:formulationAndAnalysis} the abstract problem formulation is presented along with an analysis of the conditioning of unfitted finite element methods. This section also briefly reviews the concepts underlying the SIPIC preconditioner developed in \cite{SIPIC}, which are essential for extending the preconditioner to non-SPD and mixed systems.  The CbAS preconditioning technique we propose here is described in Section~\ref{sec:preconIntro}. In Section~\ref{sec:results} we show the numerical results of the preconditioner, and concluding remarks are given in Section~\ref{sec:conclusion}.

\section{Problem formulations and conditioning analysis}
In Section~\ref{sec:formulation} we introduce the definitions and formulations that will be used throughout this manuscript, after which we present an analysis of the conditioning of unfitted finite element methods in Section~\ref{sec:conditioning}. The essential aspects of the algebraic SIPIC preconditioning scheme developed in \cite{SIPIC}, including its main restrictions, are reviewed in Section~\ref{sec:recap}.
\label{sec:formulationAndAnalysis}
\subsection{Formulations and definitions}
\begin{figure}
  \begin{center}
  \includegraphics[width=0.65\textwidth]{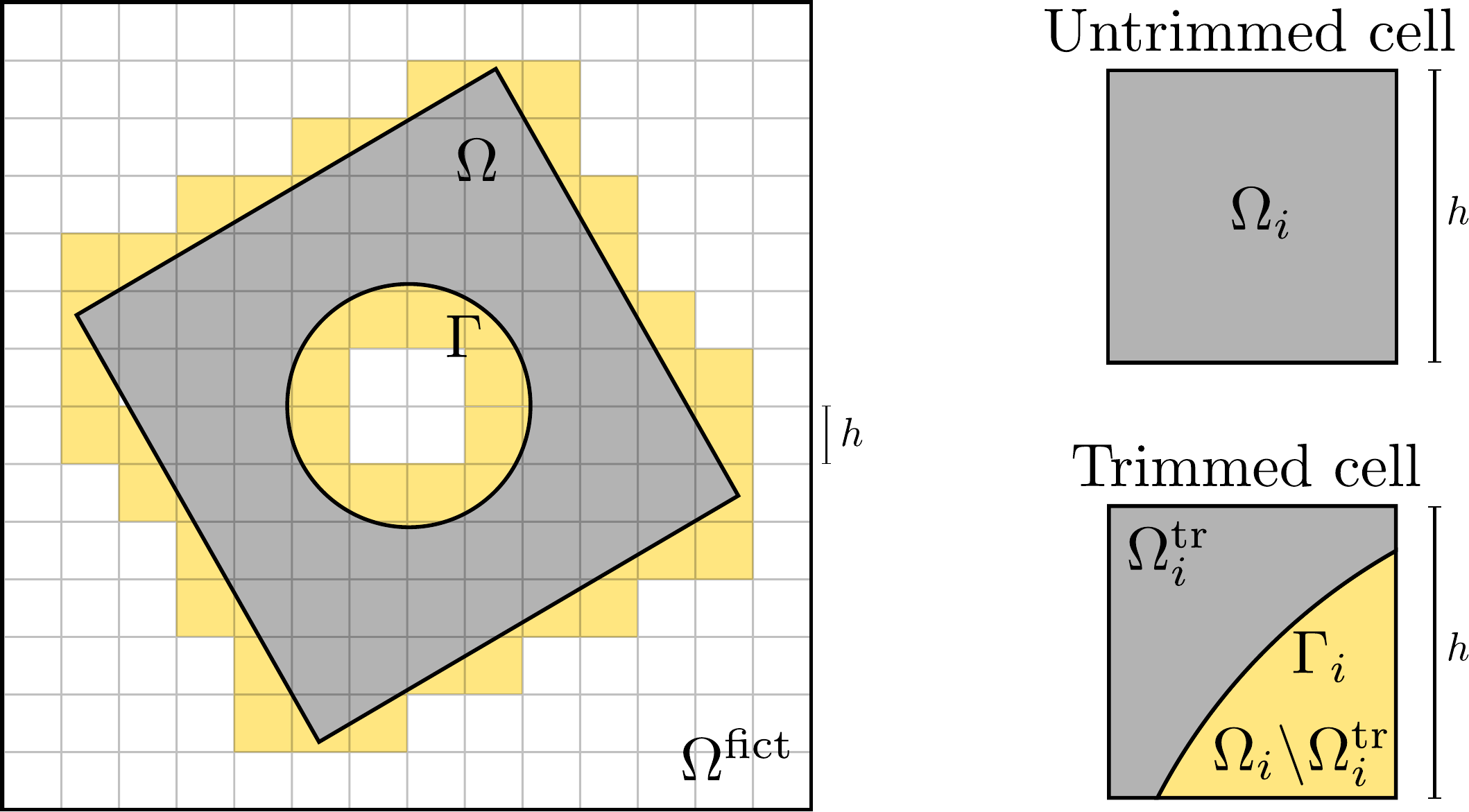}
  \caption{A geometrically complex domain $\Omega$ that is encapsulated by a fictitious domain $\Omega_{\rm fict}$ resulting in a rectilinear embedding domain $\Omega\cup\Omega_{\rm fict}$ that is simple to discretize by a tensor product grid.\label{fig:domain}}
  \end{center}
\end{figure}

We consider single-field and mixed partial differential equations posed on a $d$-dimensional domain $\Omega \subset \mathbb{R}^d$ $\left(d\in\{2,3\}\right)$, which we refer to as the physical domain. The boundary $\Gamma = \partial \Omega$ is partitioned in the complementary parts $\Gamma^D$ and $\Gamma^N$ ($\Gamma^D \cup \Gamma^N = \Gamma$ and $\Gamma^D \cap \Gamma^N = \emptyset$) for the imposition of Dirichlet and Neumann boundary conditions, respectively. The physical domain is encapsulated by a geometrically simple domain $\Omega \cup \Omega_{\rm fict}$, on which a tensor product grid can be defined (Figure~\ref{fig:domain}). Herein we consider uniform grids, but the condition number analysis and proposed preconditioning technique are not restricted to this. On the encapsulating domain $\Omega \cup \Omega_{\rm fict}$ we define multivariate isogeometric B-spline bases $\mathcal{S}_\alpha^p(\Omega \cup \Omega_{\rm fict})$, which are tensor products of univariate B-spline bases of degree $p$ and regularity $\alpha$ \cite{Hughes2005,Cottrell2009}. Note that, in principle, both $p$ and $\alpha$ can be different per spatial direction and can even vary locally. In this work we restrict ourselves to the same global $p$ and $\alpha$ for every spatial direction for simplicity however. The restriction of $\mathcal{S}_\alpha^p(\Omega \cup \Omega_{\rm fict})$ to basis functions supported on $\Omega$ is denoted by  the basis $\mathcal{S}_\alpha^p(\Omega)$. The span of basis $\mathcal{S}_\alpha^p(\Omega)$ forms the finite dimensional function space $\mathcal{V}_h(\Omega)$ in which we approximate the solution to the single-field problems. The solution to two-field mixed problems is approximated in $\mathcal{V}_h(\Omega)\times\mathcal{Q}_h(\Omega)$, which comprises separate bases for the function spaces $\mathcal{V}_h(\Omega)$ and $\mathcal{Q}_h(\Omega)$. These spaces are generally naturally equipped with inner products, corresponding to the problem under consideration.

The variational form of the single-field problems we consider herein is written as:
\begin{equation}\left\{ \begin{array}{r l}
\mbox{Find } u_h \in \mathcal{V}_h(\Omega) \mbox{ such that:} & \\[0.1em]
a(v_h,u_h) = b(v_h) & \forall v_h \in \mathcal{V}_h(\Omega).
\end{array}\right.
\label{eq:weakSingleProblem}
\end{equation}
For two-field mixed problems the variational form reads:
\begin{equation}\left\{ \begin{array}{r l}
\mbox{Find } (u_h,p_h) \in \mathcal{V}_h(\Omega)\times\mathcal{Q}_h(\Omega) \mbox{ such that:} & \\[0.1em]
a_{vu}(v_h,u_h)+a_{vp}(v_h,p_h) = b_v(v_h) & \forall v_h \in \mathcal{V}_h(\Omega), \\[0.1em]
                          a_{qu}(q_h,u_h) = b_q(q_h) & \forall q_h \in \mathcal{Q}_h(\Omega).
\end{array}\right.
\label{eq:weakMixedProblem}
\end{equation}
In these formulations $a(\cdot,\cdot)$ and $b(\cdot)$ are bounded bilinear and linear forms on $\mathcal{V}_h(\Omega)$ and $\mathcal{V}_h(\Omega)\times\mathcal{Q}_h(\Omega)$, respectively. In the examples in Section~\ref{sec:results} we have $\mathcal{V}_h(\Omega)\subset H^1(\Omega)$ and $\mathcal{Q}_h(\Omega)\subset L^2(\Omega)$. We consider variational problems that are out of the scope of \cite{SIPIC}, in the sense that the discretization of \eqref{eq:weakSingleProblem} does not result in a symmetric positive definite (SPD) matrix and that the discretization of \eqref{eq:weakMixedProblem} results in a, possibly nonsymmetric, mixed and therefore indefinite matrix.

In order to compute the solution to \eqref{eq:weakSingleProblem} we define the vector $\boldsymbol{\Phi}$ consisting of all basis functions in $\mathcal{S}_\alpha^p(\Omega)$ and spanning $\mathcal{V}_h$.  For every function $v_h \in \mathcal{V}_h(\Omega)$ there exists a unique vector $\mathbf{y}\in\mathbb{R}^n$ such that $v_h = \boldsymbol{\Phi}^T\mathbf{y}$. We can therefore condense \eqref{eq:weakSingleProblem} into the linear algebraic system:
\begin{equation}
 \mathbf{A}\mathbf{x} = \mathbf{b},
 \label{eq:matrixSingle}
\end{equation}
in which $\mathbf{A} = a(\boldsymbol{\Phi},\boldsymbol{\Phi}^T)$, $\mathbf{b} = b(\boldsymbol{\Phi})$ and $u_h = \boldsymbol{\Phi}^T\mathbf{x}$. Similarly we define $\boldsymbol{\Phi}$ spanning $\mathcal{V}_h$ and $\boldsymbol{\Psi}$ spanning $\mathcal{Q}_h$ to represent the solution to \eqref{eq:weakMixedProblem}. This yields the linear system:
\begin{equation}
 \underbrace{ \begin{bmatrix} \mathbf{A}_{vu} & \mathbf{A}_{vp} \\ \mathbf{A}_{qu} & \mathbf{0} \end{bmatrix} }_{ \mathbf{A} } \underbrace{ \begin{pmatrix} \mathbf{x}_u \\ \mathbf{x}_p \end{pmatrix} }_{ \mathbf{x} } = \underbrace{ \begin{pmatrix} \mathbf{b}_v \\ \boldsymbol{b}_q \end{pmatrix} }_{ \mathbf{b} },
 \label{eq:matrixMixed}
\end{equation}
in which $\mathbf{A}_{vu} = a_{vu}(\boldsymbol{\Phi},\boldsymbol{\Phi}^T)$, $\mathbf{A}_{vp} = a_{vp}(\boldsymbol{\Phi},\boldsymbol{\Psi}^T)$, $\mathbf{A}_{qu} = a_{qu}(\boldsymbol{\Psi},\boldsymbol{\Phi}^T)$, $\mathbf{b}_v = b_v(\boldsymbol{\Phi})$, $\mathbf{b}_q = b_q(\boldsymbol{\Psi})$, $u_h = \boldsymbol{\Phi}^T\mathbf{x}_u$ and $p_h = \boldsymbol{\Psi}^T\mathbf{x}_p$.
\label{sec:formulation}
\subsection{Conditioning of immersed finite element methods}
We briefly review the definitions and properties related to conditioning and iterative solvers that are essential for the ensuing developments. The reader is referred to, \emph{e.g.,}\ \cite{Saad} and references therein for more detailed information.

When solving a linear system of the form \eqref{eq:matrixSingle} or \eqref{eq:matrixMixed}, the condition number is an important property of the matrix $\mathbf{A}$. The Euclidean (2-norm) condition number is defined as:
\begin{equation}
 \kappa_2(\mathbf{A}) = \|\mathbf{A}\|_2 \|\mathbf{A}^{-1}\|_2 \geq 1, 
\label{eq:conditionDefinition}\end{equation}
with induced Euclidean norm:
\begin{equation}
 \|\mathbf{A}\|_2 = \max_{\mathbf{x}} \frac{\|\mathbf{A}\mathbf{x}\|_2}{\|\mathbf{x}\|_2}.
\label{eq:normDefinition}\end{equation}
Small condition numbers indicate a well-conditioned system. Iterative solvers generally display more efficient convergence behavior for well-conditioned systems. In iterative solvers, the (unknown) approximation error $\|\mathbf{x} - \mathbf{x}_i\|_2$ is also better estimated by the (known) residual $\|\mathbf{b} - \mathbf{A} \mathbf{x}_i\|_2$ for well-conditioned systems. SPD systems are frequently solved by the Conjugate Gradient (CG) method. The convergence of this method directly depends on the condition number. For systems that are not SPD, the iterative method of choice is the Generalized Minimal RESidual (GMRES) method \cite{Saad}. The convergence of GMRES does not just depend on the condition number, but a more direct dependence exists on the positioning of the eigenvalues in the complex plane. As a measure of performance of GMRES we consider the eigenvalue ratio:
\begin{equation}
 \rho(\mathbf{A}) = \frac{|\lambda|_{\rm max}}{|\lambda|_{\rm min}} \geq 1,
\end{equation}
with $|\lambda|_{\rm max}$ and $|\lambda|_{\rm min}$ the eigenvalues with the largest and smallest magnitudes. It should be noted that for symmetric matrices $\kappa_2(\mathbf{A}) = \rho(\mathbf{A})$ and that for nonsymmetric matrices $\kappa_2(\mathbf{A}) \geq \rho(\mathbf{A})$.

Immersed finite element methods are known to yield severely ill-conditioned systems when the system contains elements that only intersect the physical domain $\Omega$ on a small fraction of their volume, \emph{i.e.,}\ cut elements with small volume fractions. The volume fraction $\eta_i$ of element $\Omega_i$ (Figure~\ref{fig:domain}) is defined as the fraction of the element that intersects the physical domain $\Omega$:
\begin{equation}
\eta_i = \frac{|\Omega_i\cap\Omega|}{|\Omega_i|} = \frac{|\Omega_i^{\rm tr}|}{|\Omega_i|}.
\end{equation}
The smallest volume fraction in the system is denoted by $\eta = \min_i \eta_i$.

From \eqref{eq:conditionDefinition} and \eqref{eq:normDefinition} it follows that the condition number depends on the following maxima:
\begin{equation}
 \|\mathbf{A}\|_2 = \max_{\|\mathbf{y}\|_2=1} \| \mathbf{A}\mathbf{y} \|_2 =  \max_{\|\mathbf{y}\|_2=1} \| a(\boldsymbol{\Phi},\boldsymbol{\Phi}^T\mathbf{y}) \|_2 = \max_{\substack{ v_h=\boldsymbol{\Phi}^T\mathbf{y}, \\ \|\mathbf{y}\|_2=1}} \| a(\boldsymbol{\Phi},v_h) \|_2,
 \label{eq:normmax}
\end{equation}
\begin{equation}
 \|\mathbf{A}^{-1}\|_2 = \max_{\|\mathbf{y}\|_2=1} \frac{1}{\| \mathbf{A}\mathbf{y} \|_2} = \Bigg( \min \limits_{\substack{ v_h=\boldsymbol{\Phi}^T\mathbf{y},\\\|\mathbf{y}\|_2=1}} \| a(\boldsymbol{\Phi},v_h) \|_2 \Bigg)^{-1},
 \label{eq:normmin}
\end{equation}
with $\mathbf{A}=a(\boldsymbol{\Phi},\boldsymbol{\Phi}^T)$. The norm and (in magnitude) largest eigenvalue of $\mathbf{A}$ are generally unaffected by the trimmed elements, \emph{viz.}\ these are of approximately the same magnitude for fitted (mesh-conforming) and unfitted (immersed) methods. For most variational problems the maximum in \eqref{eq:normmax} is obtained by an oscillatory function with the highest frequency that can be captured by the grid, such as proven for the Poisson problem in \cite{Johnson}. The norm and in magnitude largest eigenvalue of $\mathbf{A}^{-1}$ on the other hand behave very differently for mesh-fitting and immersed methods. Note that the largest eigenvalue of $\mathbf{A}^{-1}$ is equal to the reciprocal smallest eigenvalue of $\mathbf{A}$. For most variational problems, mesh-fitting methods with regular grids obtain the maximum in \eqref{eq:normmin} by a function that closely resembles the eigenfunction of the operator with the lowest eigenfrequency, see \emph{e.g.,}\ \cite{Johnson} for the Poisson problem. For immersed methods on the other hand, the maximum in \eqref{eq:normmin} is generally attained by a function that is only supported on the element with the smallest volume fraction. This directly follows from \eqref{eq:normmin}: as $a(\boldsymbol{\Phi},v_h)$ is a linear operator on $v_h$, $\| a(\boldsymbol{\Phi},v_h) \|_2$ will be small for a small function $v_h$. When the system contains a cut element with a very small volume fraction, a function $v_h$ that is only supported on that element will be small in any norm. The volume fraction of this element does not affect the norm of the corresponding vector $\|\mathbf{y}\|$, however. Therefore there will be functions in the system with $\|v_h\| \ll \|\mathbf{y}\|$ in any suitable norms. This yields very large values for $\|\mathbf{A}^{-1}\|$ which causes the ill-conditioning of immersed methods.

With high-order bases, the effect is exacerbated. As shown in \cite[Section~3.2]{SIPIC}, polynomial bases generally contain a function of the form:
\begin{equation}
v_h = \prod_{j=1}^d \left(\frac{x_j-\hat{x}_j}{h}\right)^{p},
\label{eq:smallfunc}
\end{equation}
that is only supported on the element with the smallest volume fraction and that corresponds to a coefficient vector with a magnitude of order one. In \eqref{eq:smallfunc}, $h$ is the mesh size, $p$ the order of the discretization, and $\hat{x}$ a point on the boundary of the element. For small volume fractions $|x - \hat{x}|_e \ll h$ such that $v_h \ll 1$. Note that we assume $h$ and $p$ to be isotropic. However, a similar relation holds for anisotropic $h_j$ and $p_j$ per spatial direction. It follow from \eqref{eq:smallfunc} that high-order bases contain even smaller functions than lower order bases (both with a corresponding vector of norm $\|\mathbf{y}\|=1$). For symmetric and elliptic second order problems, we have used this analysis to derive a lower bound for the condition number by the scaling relation \cite[Section~3.2]{SIPIC}:
\begin{equation}
 \kappa_2(\mathbf{A}) \geq C \eta^{(2p+1-2/d)},
 \label{eq:spdscalingrelation}
\end{equation}
for some constant $C>0$ independent of $\eta$. 

In the derivation of \eqref{eq:spdscalingrelation} a uniform grid was assumed, but a similar derivation can be established for nonuniform grids. Similar relations can also be derived for different single-field and mixed partial differential equations. The rate generally depends on the equation under consideration, and for mixed equations even the specific pair of function spaces (\emph{e.g.,}\ Taylor-Hood, Raviart-Thomas, \emph{etc.}) has an effect. The scaling factor $2p$ in this rate appears to be generic however. The derivation of the exact scaling relations for all variational problems considered in Section~\ref{sec:results} is beyond the scope of this work, but the existence of such scaling relations is supported by the numerical results.
\label{sec:conditioning}
\subsection{Algebraic preconditioning and limitations}
If a system $\mathbf{A}\mathbf{x}=\mathbf{b}$ is preconditioned, one instead solves one of the the equivalent systems: 
\begin{equation}\begin{array}{r l}
\mathbf{S}\mathbf{A} \mathbf{x} = \mathbf{S}\mathbf{b} \hspace{.15cm}\rule{0pt}{0pt} \hspace{.3cm}\rule{0pt}{0pt} & \mbox{(left preconditioning),} \vspace{.2cm} \\
\left.\begin{aligned} \mathbf{A}\mathbf{S} \mathbf{y} & = \mathbf{b} \\ \mathbf{x} & = \mathbf{S} \mathbf{y} \hspace{.15cm}\rule{0pt}{0pt} \end{aligned}\right\} & \mbox{(right preconditioning),} \vspace{.2cm} \\
\left.\begin{aligned} \mathbf{S}^{\textsc{l}}\mathbf{A}\mathbf{S}^{\textsc{r}} \mathbf{y} & = \mathbf{S}^{\textsc{l}}\mathbf{b} \\ \mathbf{x} & = \mathbf{S}^{\textsc{r}} \mathbf{y} \end{aligned}\right\} & \mbox{(left and right preconditioning).}
\end{array}\end{equation}
The matrices $\mathbf{S}$, $\mathbf{S}^{\textsc{l}}$ and $\mathbf{S}^{\textsc{r}}$ are chosen such that the preconditioned matrices $\mathbf{S}\mathbf{A}$, $\mathbf{A}\mathbf{S}$ or $\mathbf{S}^{\textsc{l}}\mathbf{A}\mathbf{S}^{\textsc{r}}$ have a smaller condition number or eigenvalue ratio than matrix $\mathbf{A}$ itself. In the examples in Section~\ref{sec:results} we restrict ourselves to left preconditioning, as left and right preconditioning are equivalent with respect to eigenvalues and as also $\mathbf{S}^{\textsc{l}}\mathbf{A}\mathbf{S}^{\textsc{r}}$ has the same eigenvalues as $\mathbf{S}^{\textsc{r}}\mathbf{S}^{\textsc{l}}\mathbf{A}$. Generally $\mathbf{S}$ or $\mathbf{S}^{\textsc{r}}\mathbf{S}^{\textsc{l}}$ is a sparse approximation of $\mathbf{A}^{-1}$, as this implies $\kappa_2(\mathbf{S}\mathbf{A}) \approx \kappa_2(\mathbf{I}) = 1$ and $\rho(\mathbf{S}\mathbf{A}) \approx \rho(\mathbf{I}) = 1$.

The underlying mechanism of the ill-conditioning of immersed methods are functions $\boldsymbol{\Phi}^T\mathbf{y} = v_h \in \mathcal{V}_h(\Omega)$ that are only supported on cut elements with small volume fractions such that $\| v_h \| \ll \| \mathbf{y} \|$. The conceptual idea of the algebraic SIPIC preconditioner is to preclude the aforementioned effect by preconditioning the basis $\boldsymbol{\Phi}$. This is done by instigating a change of basis with a nonsingular preconditioning matrix $\mathbf{S}^\frac12$:
\begin{equation}
 \overline{\boldsymbol{\Phi}} = \mathbf{S}^\frac12 \boldsymbol{\Phi}.
\end{equation}
Basis $\overline{\boldsymbol{\Phi}}$ spans the same space as $\boldsymbol{\Phi}$ due to the nonsingularity of $\mathbf{S}^\frac12$, but by adequate construction of this preconditioner we prevent the occurrence of functions for which the magnitude of the corresponding coefficient vector is much larger than the function itself. For the system matrix this implies:
\begin{equation}
 \overline{\mathbf{A}} = a(\overline{\boldsymbol{\Phi}},\overline{\boldsymbol{\Phi}}^T) = a(\mathbf{S}^\frac12\boldsymbol{\Phi},\boldsymbol{\Phi}^T\mathbf{S}^{\frac12,T}) = \mathbf{S}^\frac12a(\boldsymbol{\Phi},\boldsymbol{\Phi}^T)\mathbf{S}^{\frac12,T} = \mathbf{S}^\frac12\mathbf{A}\mathbf{S}^{\frac12,T},
\end{equation}
which has the same eigenvalues as the left preconditioned system:
\begin{equation}
 \mathbf{S}\mathbf{A} \mathbf{x} = \mathbf{S}\mathbf{b},
\end{equation}
with:
\begin{equation}
 \mathbf{S} = \mathbf{S}^{\frac12,T}\mathbf{S}^\frac12.
\end{equation}
A simple choice for $\mathbf{S}^\frac12$ (and consequently $\mathbf{S}$) is a diagonal matrix that linearly scales the basis functions. For positive definite systems this scaling can be derived from the square root of the diagonal of $\mathbf{A}$, such that all basis functions in $\overline{\boldsymbol{\Phi}}$ are equal in the operator norm (\emph{i.e.,}\ $\| \phi \|_a^2 = a( \phi,\phi)=1~\forall \phi \in \overline{\boldsymbol{\Phi}}$). This preconditioner is also known as diagonal scaling or Jacobi preconditioning.

In \cite[Section 4.1.2]{SIPIC} it has been shown that linear scaling is generally not sufficient to obtain proper conditioning. On elements with small volume fractions, basis functions do not only become very small, but also higher order terms become less significant. As a result, scaled basis functions can be very similar to each other and therefore become almost linearly dependent. To bypass this effect, the SIPIC preconditioning scheme algebraically identifies and orthonormalizes the functions that are almost linearly dependent. This results in adding a few off-diagonal terms, leading to a very sparse (almost diagonal) preconditioner. After a permutation this leads to a block lower diagonal matrix of the form displayed in Figure~\ref{fig:blockLower}.

The scaling and orthonormalization operations in the construction of the SIPIC preconditioner do require the variational form to be an inner product, which restricts this preconditioning approach to SPD problems. This precludes application of SIPIC to most flow problems. Furthermore, even though the method to identify the almost linearly dependent functions has shown to be very effective for smooth B-spline bases, this method is not always adequate for high-order bases with low regularity. The reason for this is that in such bases it can occur that sets of more than two basis functions are almost linearly dependent, but that in the SIPIC construction based on all possible pairs of functions in this set this dependence is not identified. This results in potential suboptimal performance of the SIPIC preconditioner for such systems. The CbAS preconditioner proposed in Section~\ref{sec:preconIntro} does not suffer from these restrictions, and can hence be conceived of as a generalization and improvement of SIPIC.
\label{sec:recap}

\section{Preconditioning immersed finite element methods}
In this section we introduce the \emph{Connectivity-based Additive-Schwarz} (CbAS) preconditioner. CbAS enables preconditioning of immersed finite element approximations of flow problems by generalization of two aspects of SIPIC. In Section~\ref{sec:preconSingle} we first propose a systematic way of constructing and preconditioning connectivity-based matrix blocks. This makes CbAS applicable to high-order bases with low regularity and to problems with non-SPD matrices. In Section~\ref{sec:preconCost} we comment on the computational cost of the CbAS preconditioner. In Section~\ref{sec:preconMixed} we propose a procedure to optimally apply CbAS to mixed problems.\label{sec:preconIntro}
\subsection{The CbAS preconditioner for single-field problems}
\begin{figure}

  \begin{subfigure}[t]{.49\textwidth}
    \centering
    \includegraphics[height=4cm]{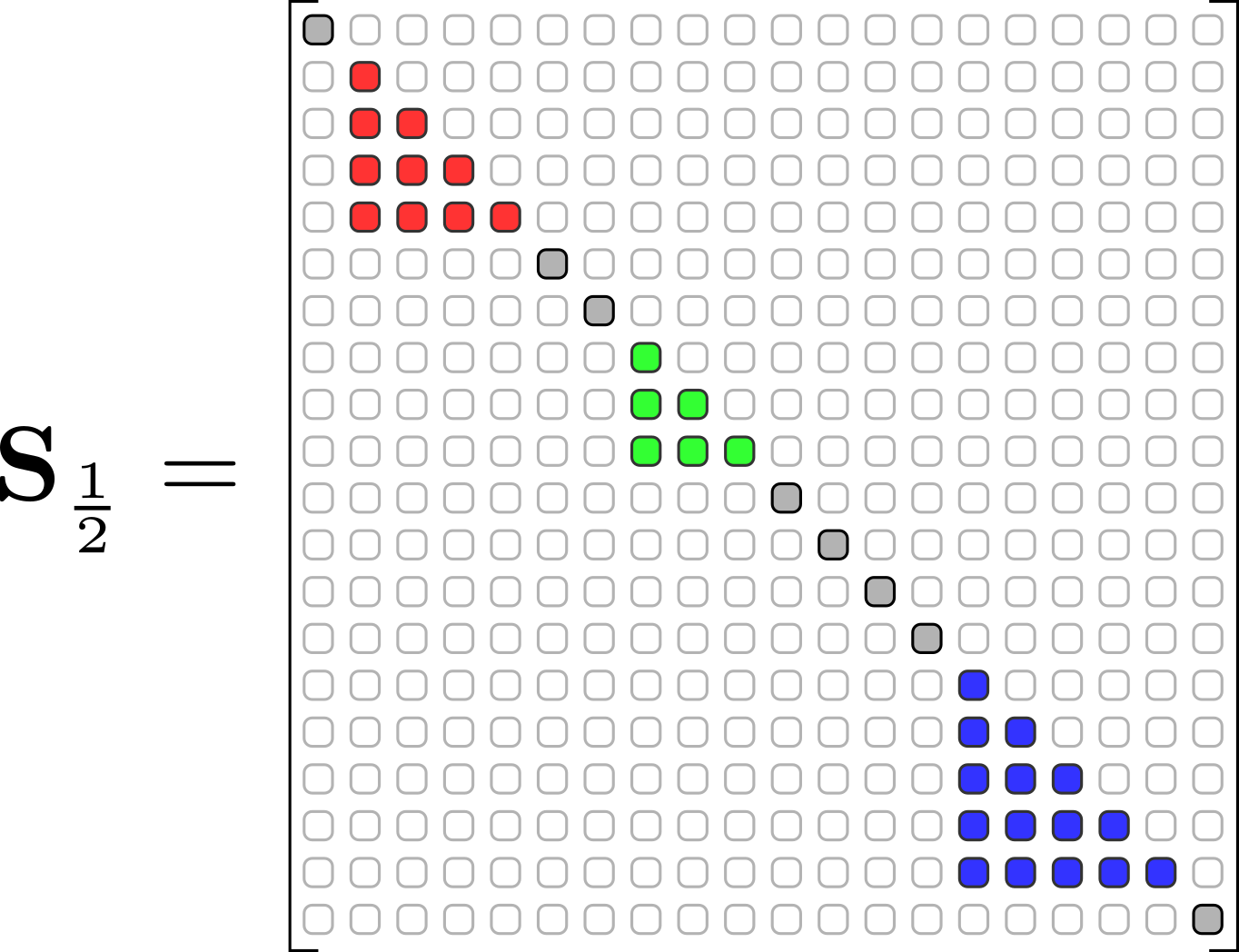}
    \caption{}\label{fig:blockLower}
  \end{subfigure}
  \begin{subfigure}[t]{.49\textwidth}
    \centering
    \includegraphics[height=4cm]{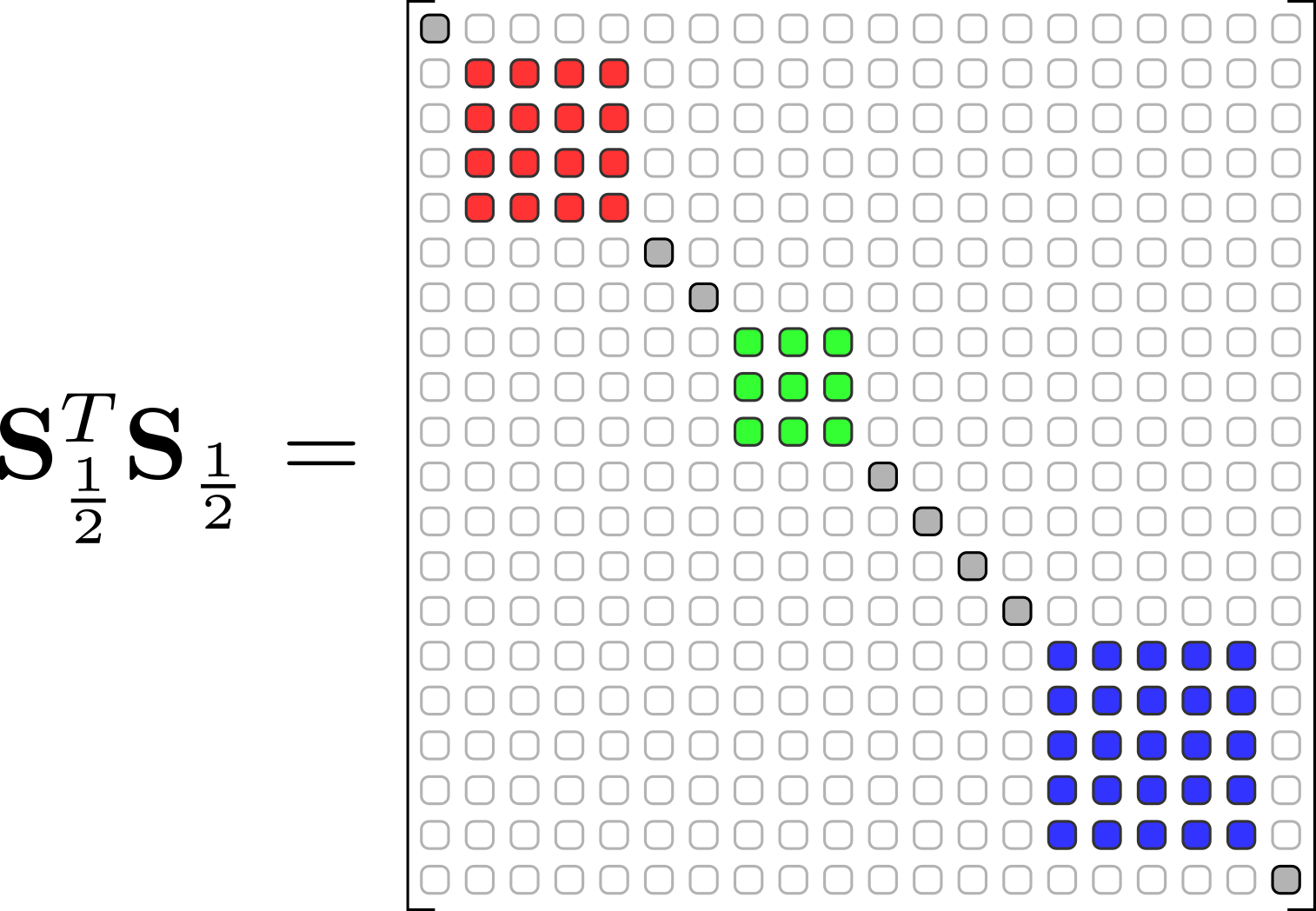}
    \caption{}\label{fig:blockDistinct}
  \end{subfigure}

\caption{Illustrative sparsity pattern of (a) the lower diagonal preconditioner of the basis $\mathbf{S}^\frac12$, and (b) the block-diagonal preconditioner $\mathbf{S} = \mathbf{S}^{\frac12,T}\mathbf{S}^\frac12$. Note that both matrices are permuted and separable as they contain non-overlapping blocks.}
\end{figure}

To motivate the CbAS preconditioner, we first reconsider some of the properties underlying the construction of the SIPIC preconditioner. Similar to $\mathbf{S}^\frac12$ in Figure~\ref{fig:blockLower}, the left SIPIC preconditioner $\mathbf{S} = \mathbf{S}^{\frac12,T}\mathbf{S}^\frac12$ can be permuted to a separable block diagonal matrix, see Figure~\ref{fig:blockDistinct}. Furthermore, $\mathbf{S}$ is equal to the inverse of the matrix that would result from restricting the SPD matrix $\mathbf{A}$ to the same separable blocks of nonzero entries as $\mathbf{S}$, \emph{viz.}\ the blocks of $\mathbf{S}$ are the inverses of the blocks of $\mathbf{A}$. This is because $\mathbf{S}^\frac12$ corresponds to the inverse Cholesky decomposition of this restriction of $\mathbf{A}$, \emph{viz.}\ the lower diagonal blocks of $\mathbf{S}^\frac12$ are the inverse Cholesky decompositions of the blocks of $\mathbf{A}$. Computing these block-wise Cholesky decompositions by the Gram-Schmidt orthonormalization procedure is one of the aspects that restrict the SIPIC preconditioning technique to SPD matrices. If we omit computing $\mathbf{S}^\frac12$ but directly compute the blocks of $\mathbf{S}$ by inverting the blocks of $\mathbf{A}$, this restriction to SPD matrices is relaxed to $\mathbf{A}$ being (block-wise) nonsingular. This condition is generally satisfied for well-posed variational forms. Therefore -- except for identifying the functions that are almost linearly dependent -- this generalizes the SIPIC preconditioning scheme to non-SPD matrices.

\begin{figure}
  \begin{center}
  \includegraphics[height=4cm]{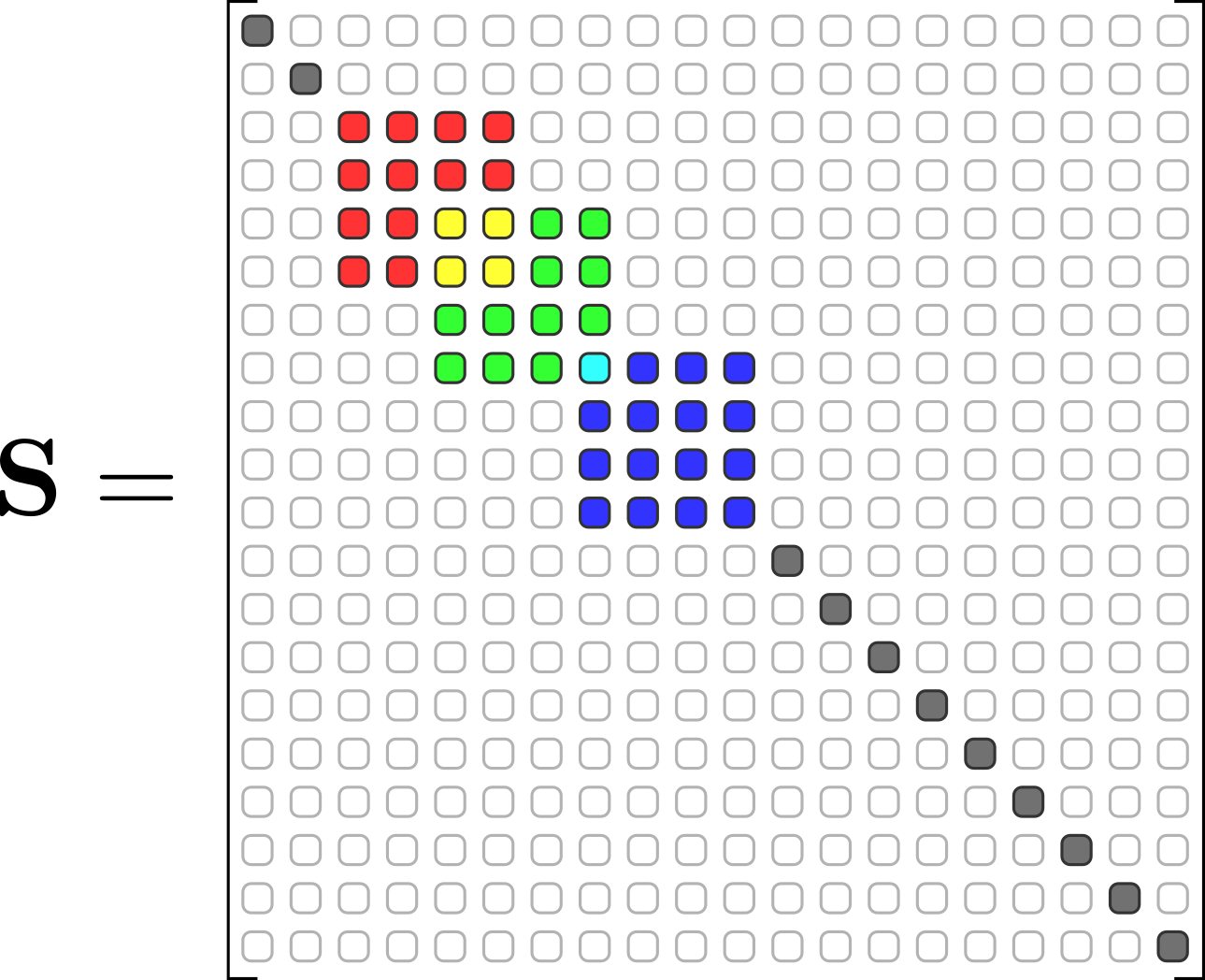}
  \caption{Illustrative sparsity pattern of a block-diagonal Additive-Schwarz preconditioner $\mathbf{S}$ constructed by \eqref{eq:ASconstruct}. Values of the red and green blocks and green and blue blocks are summed at the yellow and respectively cyan indices. Note that because of the overlap these blocks are not separable anymore in the assembled preconditioning matrix.\label{fig:blockOverlap}}
  \end{center}
\end{figure}

The procedure in \cite[Section~4.1.2]{SIPIC} to identify the near linear dependencies in SIPIC is also restricted to SPD matrices. Additionally, the procedure has shown to be less robust for high-order bases with low regularity than for the smooth B-spline bases considered in \cite{SIPIC}. A more robust procedure results, however, if the identification of these dependencies is approached in a similar manner as in domain decomposition preconditioning (see \emph{e.g.,}\ \cite{Smith1996,Toselli2005}). It is well known that in the framework of Additive-Schwarz preconditioners, it is not necessary that these blocks are separable (\emph{i.e.,}\ non-overlapping), as is the case for SIPIC (see Figures~\ref{fig:blockDistinct} and \ref{fig:blockOverlap}). There exists a rich literature on Additive-Schwarz preconditioners for finite element methods (see \emph{e.g.,}\ \cite{Ferencz1998}) and recent work for isogeometric analysis \cite{BeiraoDaVeiga2012,BeiraoDaVeiga2013Elasticity,BeiraoDaVeiga2014}.

To provide a basis for the CbAS preconditioner developed herein we first recall some elementary aspects of Additive-Schwarz preconditioning. We denote the size of matrix $\mathbf{A}$ by $n$, the number of blocks by $I$, the index of a block by $i \leq I$ and the rank of block $i$ by $m_i$. For every block a projection matrix $\mathbf{P}_i \in \mathbb{R}^{n \times m_i}$ is constructed, consisting of the unit vectors of the indices of the functions in block $i$. For example, if block $i$ contains the functions with indices $\{\alpha,\beta,\gamma\}$, then $\mathbf{P}_i (x_\alpha, x_\beta, x_\gamma)^T$ creates a vector of length $n$ with $\{x_\alpha,x_\beta,x_\gamma\}$ at the indices $\{\alpha,\beta,\gamma\}$ as the only nonzero entries. The transpose of projection matrix $\mathbf{P}_i$ is a reduction matrix, \emph{i.e.,}\ $\mathbf{P}_i^T\mathbf{y} = (y_\alpha,y_\beta,y_\gamma)^T$ for any vector $\mathbf{y}$ of length $n$. The restriction of matrix $\mathbf{A}$ to block $i$ is denoted by $\mathbf{A}_i = \mathbf{P}_i^T\mathbf{A}\mathbf{P}_i$. Assuming that the blocks are invertible -- \emph{e.g.,}\ $\mathbf{A}$ derives from a coercive bilinear form -- if the blocks $\mathbf{A}_i$ are sufficiently small then it is feasible to invert each $\mathbf{A}_i$ (directly or approximately) to form $\mathbf{A}_i^{-1}$. The inverse of $\mathbf{A}_i$ is projected into an $n \times n$ matrix using the projection and restriction matrices, \emph{i.e.,}\ $\mathbf{P}_i\mathbf{A}_i^{-1}\mathbf{P}_i^T$. This results in a sparse matrix with only $m_i^2$ nonzero entries at the cross indices of block $i$. Finally these sparse matrices are summed to form the preconditioner:
\begin{equation}
 \mathbf{S} = \sum_{i=1}^I \mathbf{P}_i\underbrace{\left(\mathbf{P}_i^T\mathbf{A}\mathbf{P}_i\right)^{-1}}_{\mathbf{A}_i^{-1}}\mathbf{P}_i^T.
 \label{eq:ASconstruct}
\end{equation}
Because $\mathbf{S}$ is a summation of block-wise inverses of $\mathbf{A}$, it is intuitive to interpret $\mathbf{S}$ as a sparse approximation of $\mathbf{A}^{-1}$. Every index must be in at least one of the blocks, as otherwise $\mathbf{S}$ contains empty rows and columns and therefore exhibits a null space. If this requirement is satisfied and all the blocks are non-overlapping, the preconditioner $\mathbf{S}$ is again equal to the inverse of the restriction of $\mathbf{A}$ to the same blocks. For SPD matrices the procedure in \eqref{eq:ASconstruct} then yields the same preconditioner as the local orthonormalization procedure in SIPIC. The Additive-Schwarz procedure can therefore be conceived of as a generalization of the local orthonormalization procedure in SIPIC in the sense that it can also treat non-SPD matrices and overlapping (non-separable) blocks.

If $\mathbf{A}$ is symmetric positive definite the Additive-Schwarz lemma (commonly referred to as \emph{Lions' lemma}) holds \cite{Matsokin1985,Lions1988} (see \cite{Smith1996,Toselli2005} for this specific form):
\begin{equation}
 \mathbf{y}^T \mathbf{S}^{-1} \mathbf{y} = \min_{\mathbf{y} = \sum\limits_{j=1}^{I} \mathbf{P}_j \mathbf{y}_j} \hspace{.4cm} \sum_{i=1}^I \mathbf{y}_i^T \mathbf{A}_i \mathbf{y}_i = \min_{\mathbf{y} = \sum\limits_{j=1}^{I} \mathbf{P}_j \mathbf{y}_j} \hspace{.4cm} \sum_{i=1}^I \left( \mathbf{P}_i \mathbf{y}_i \right)^T \mathbf{A} \left( \mathbf{P}_i \mathbf{y}_i \right) \quad \forall \mathbf{y} \in \mathbb{R}^n.
 \label{eq:ASlemma}
\end{equation}
This lemma indicates that Rayleigh quotients with $\mathbf{S}^{-1}$ are closely related to Rayleigh quotients with $\mathbf{A}$, emphasizing the intuitive interpretation of $\mathbf{S}$ as a sparse approximation of $\mathbf{A}^{-1}$. As discussed in Section~\ref{sec:conditioning}, the underlying mechanism of ill-conditioning of immersed methods pertains to trimmed and therefore small or almost linearly dependent basis functions. Such basis functions can together form a function $v_h$ and corresponding coefficient vector $\mathbf{y}$ with $\| v_h \| \ll \| \mathbf{y} \|$, such that in the SPD case $\mathcal{F}(v_h,v_h) = \mathbf{y}^T\mathbf{A}\mathbf{y} \ll \mathbf{y}^T\mathbf{y}$. From \eqref{eq:ASlemma} one can see that if such almost linearly dependent functions are in the same block, they will also form a small Rayleigh quotient with $\mathbf{S}^{-1}$. The Additive-Schwarz preconditioner $\mathbf{S}$ then effectively targets the fundamental conditioning problem of immersed methods. The Additive-Schwarz lemma in \eqref{eq:ASlemma} is in principle restricted to SPD systems, but the procedure has been analyzed extensively for nonsymmetric and indefinite systems \cite{Cai1991,Cai1992,Widlund1992,Cai1993,Cai1994,Chan1994,Smith1996,Cai1998,Toselli2005,Sarkis2007} and has for example successfully been applied to mesh-conforming Navier-Stokes systems (see \emph{e.g.,}\ \cite{Bazilevs2010Weak}).

A fundamental aspect in the preconditioning of immersed finite element methods based on \eqref{eq:ASconstruct} is the way in which the blocks are constructed. The selection should be done in such a way that almost linearly dependent functions are in the same block. Evidently, there are various ways to satisfy this criterion for the selection procedure. Because almost linearly dependent functions necessarily have overlapping supports, the most straightforward way to do this is to use the connectivity of the basis functions and to devise one block for every element, containing all the functions that are supported on it. On elements that are not trimmed the functions are generally sufficiently orthogonal, and therefore it is more efficient to only devise blocks for trimmed elements. This requires a separate step for functions that are only supported on non-trimmed elements. For these functions we apply standard diagonal scaling, which can be conceived of as assigning a separate $1 \times 1$ block to each such function. We will refer to the corresponding Additive-Schwarz preconditioner \eqref{eq:ASconstruct} with connectivity-based blocks for all trimmed elements and diagonal scaling for functions that are not contained in any of these blocks as the Connectivity-based Additive-Schwarz (CbAS) preconditioner.

We do note that the effectiveness of devising one block for every trimmed element does not necessarily extend to hierarchically refined grids, because of the connectivity structure for such grids. Also, the CbAS preconditioner may be enhanced with more advanced preconditioning techniques such as multigrid-type approaches (see \emph{e.g.,}\ \cite{Cai1991,Cai1992,Widlund1992,Smith1996,Toselli2005}).
\label{sec:preconSingle}
\subsection{Computational cost}
The computational cost of applying CbAS is limited. The number of basis functions that are supported per element is $(p+1)^d$ for every solution variable, with $p$ the order of the discretization and $d$ the number of dimensions. The cost of inverting the local matrices to construct the preconditioner is therefore limited. For fine grids, the number of trimmed elements is much smaller than the total number of elements. Typically, the number of trimmed elements scales with $h^{-(d-1)}$ while the total number of elements scales with $h^{-d}$. The computational cost of setting up the preconditioner, storing the preconditioner, and multiplying vectors with the preconditioner is therefore only a fraction of the inevitable cost of assembling the complete system matrix, storing the matrix, and evaluating matrix-vector products. Moreover, the Additive-Schwarz procedure underlying CbAS is easily parallelizable, because the block-wise preconditioners can be computed fully separate.

It is possible to further reduce the number of blocks. We have observed similar condition numbers and eigenvalue ratios when blocks were only devised for elements with volume fractions below a certain threshold. Also these blocks can be further restricted to only the small functions (\emph{e.g.,}\ in the operator norm). However, because the computational cost involved in the preconditioning scheme described above is already small, this makes the procedure unnecessarily complicated.

\begin{figure}
  \begin{center}
  \begin{subfigure}{0.49\textwidth}
  \centering
  \includegraphics[height=4cm]{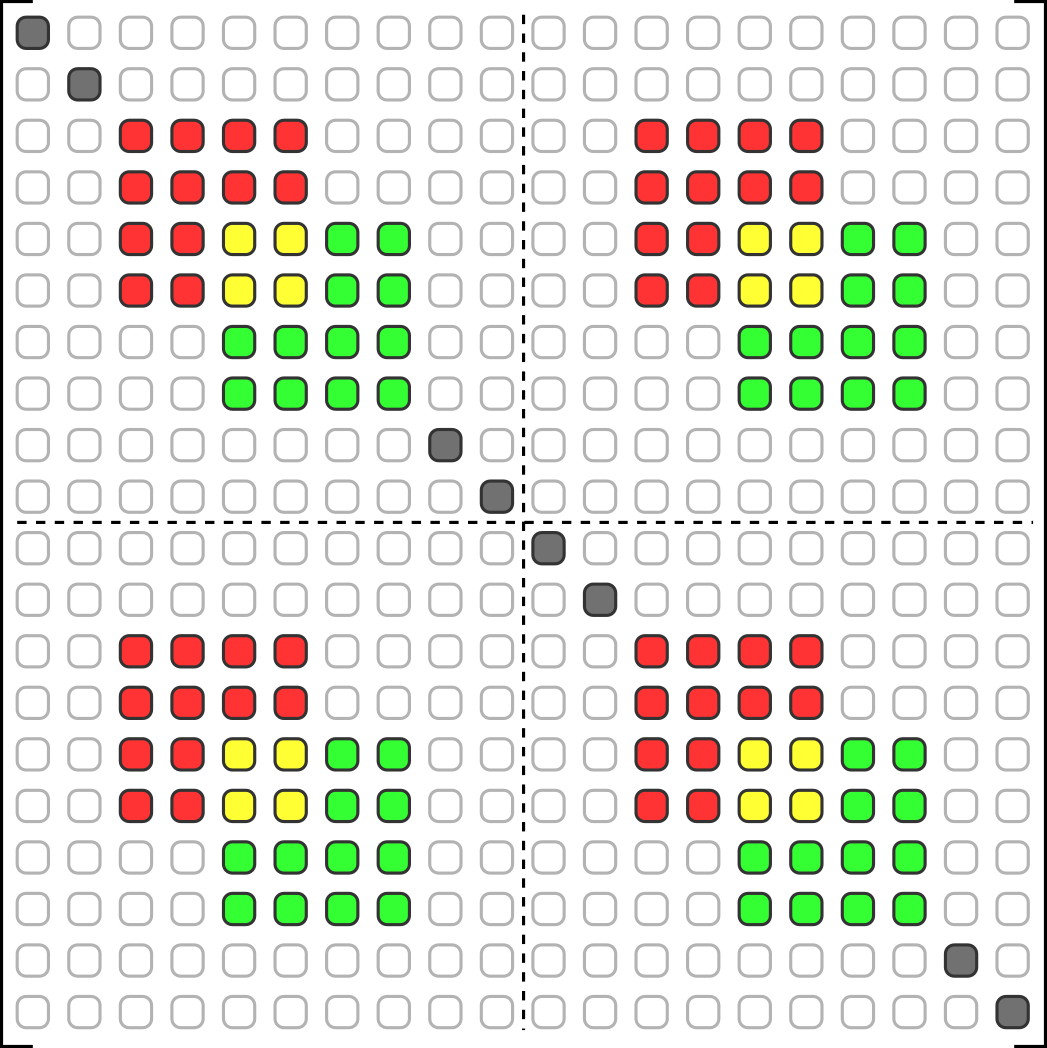}
  \caption{One block per element}
  \end{subfigure}
  \begin{subfigure}{0.49\textwidth}
  \centering
  \includegraphics[height=4cm]{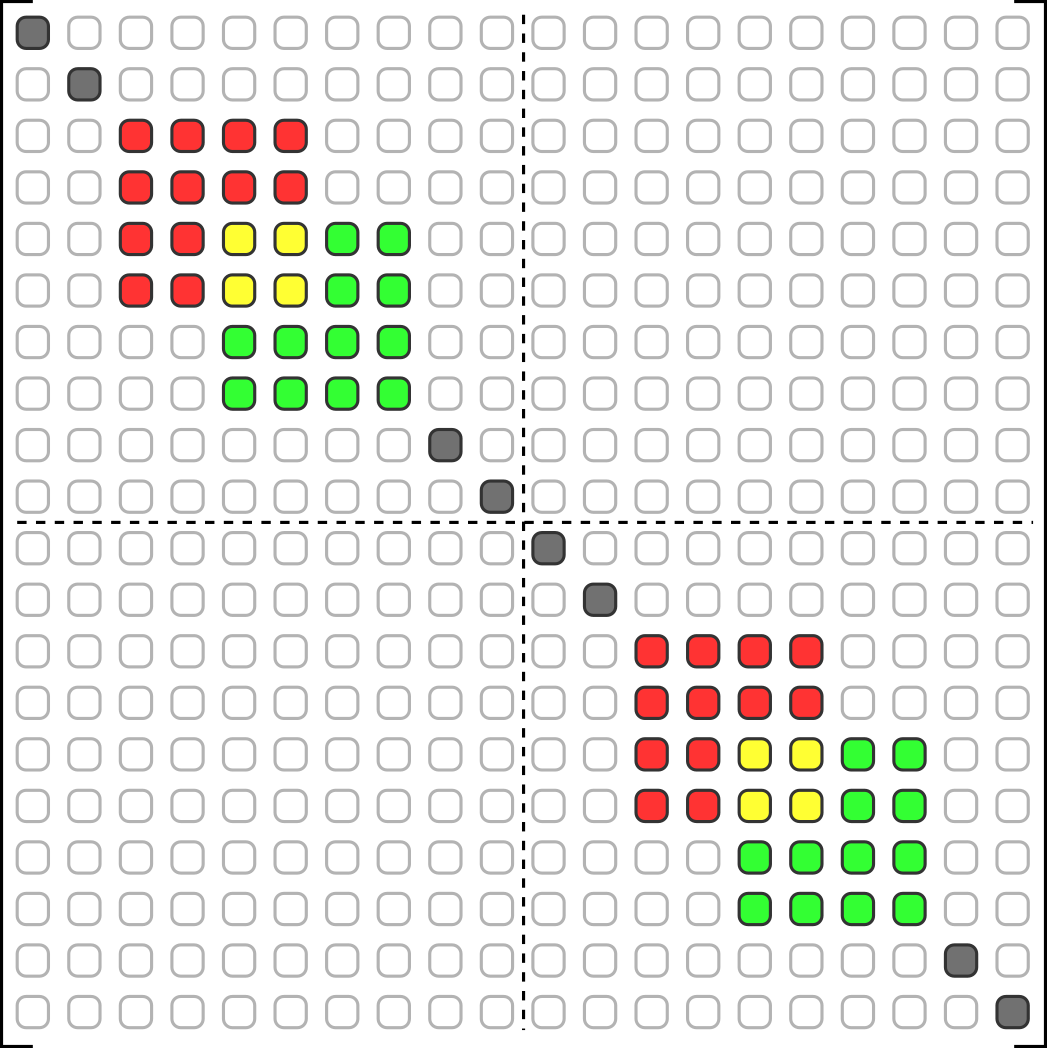}
  \caption{Separate blocks for every spatial direction}
  \end{subfigure}
  \caption{Typical sparsity patterns of the CbAS preconditioner for a two-dimensional elasticity problem.\label{fig:elasticityPattern}}
  \end{center}
\end{figure}

For multivariate systems with vector-valued solutions the computational efficiency of CbAS can be enhanced further by exploiting the structure of such problems. In the case of linear elasticity, for example, the original basis $\boldsymbol{\Phi}$ consists of separate basis functions for each of the spatial directions. Basis functions in the preconditioned basis $\overline{\boldsymbol{\Phi}} = \mathbf{S}^\frac12 \boldsymbol{\Phi}$ consist of linear combinations of basis functions of the original basis $\boldsymbol{\Phi}$. Because basis functions in separate spatial dimensions can never become linearly dependent, this does not require the preconditioned basis $\overline{\boldsymbol{\Phi}}$ to contain linear combinations of basis functions of the original basis $\boldsymbol{\Phi}$ in different spatial directions. The preconditioner $\mathbf{S}^\frac12$ and consequently $\mathbf{S}$ can therefore be restricted to block-diagonal matrices with separate blocks for every spatial direction, which is illustrated for a two-dimensional elasticity problem in Figure~\ref{fig:elasticityPattern}. This has been applied for the velocity functions in the Stokes and Navier-Stokes test cases in Section~\ref{sec:Stokes} and Section~\ref{sec:NavierStokes}.

Finally, the computational cost of the preconditioner can also be bounded when block is devised for every element, which can (by approximation) be the case when immersed finite element methods are applied to problems posed on porous domains. Because the entries in system matrix $\mathbf{A}$ are only nonzero if the functions associated with these entries have intersecting supports, the preconditioner $\mathbf{S}$ will in that case have exactly the same sparsity pattern as the system matrix. The required memory to store the preconditioner and the cost of multiplying a vector with the preconditioner will therefore be equal to the required memory and the cost of matrix-vector multiplication for the system matrix, respectively.\label{sec:preconCost}
\subsection{Extension to mixed problems}
The CbAS preconditioner described in Section~\ref{sec:preconSingle} can be extended to mixed problems, for which the solution vector represents different physical quantities. The block structure of such systems can be exploited in the construction of the preconditioner.

We restrict our considerations here to the prototypical Stokes and Navier-Stokes equations. In most variational forms of the Stokes and Navier-Stokes equations -- which involve a velocity field and a pressure field -- the pressure-pressure matrix block is empty. Hence, one cannot simply invert blocks of the pressure-pressure matrix similar to the multiple distinct spatial directions as we proposed in Section~\ref{sec:preconCost}. To develop a preconditioner for the pressure we therefore need to consider the interaction between the pressure and the velocity functions. For a symmetric variational form of the Stokes equations, symmetric preconditioning with a block-diagonal matrix with separate blocks $\mathbf{S}_u^\frac12$ and $\mathbf{S}_p^\frac12$ for the velocity and pressure, respectively, yields:
\begin{equation}
 \overline{\mathbf{A}} = \begin{bmatrix} \mathbf{S}_u^\frac12 & \mathbf{0} \\ \mathbf{0} & \mathbf{S}_p^\frac12 \end{bmatrix} \begin{bmatrix} \mathbf{A}_{vu} & \mathbf{A}_{vp} \\ \mathbf{A}_{vp}^T & \mathbf{0} \end{bmatrix} \begin{bmatrix} \mathbf{S}_u^{\frac12,T} & \mathbf{0} \\ \mathbf{0} & \mathbf{S}_p^{\frac12,T} \end{bmatrix} = \begin{bmatrix} \mathbf{S}_u^\frac12 \mathbf{A}_{vu} \mathbf{S}_u^{\frac12,T} & \mathbf{S}_u^\frac12 \mathbf{A}_{vp} \mathbf{S}_p^{\frac12,T} \\ \left( \mathbf{S}_u^\frac12 \mathbf{A}_{vp} \mathbf{S}_p^{\frac12,T} \right)^T & \mathbf{0} \end{bmatrix}.
 \label{eq:mixedpreconmatrix}
\end{equation}
The optimal preconditioner $\mathbf{S}_u$ of the velocity-velocity matrix $\mathbf{A}_{vu}$ is simply its inverse:
\begin{equation}
\mathbf{S}_u = \mathbf{S}_u^{\frac12,T} \mathbf{S}_u^\frac12 = \mathbf{A}_{vu}^{-1},
\end{equation}
and, hence:
\begin{equation}
\overline{\mathbf{A}}_{vu} = \mathbf{S}_u^\frac12 \mathbf{A}_{vu} \mathbf{S}_u^{\frac12,T} = \mathbf{I}.
\end{equation}
Considering the eigenvalue problem for the matrix in \eqref{eq:mixedpreconmatrix} with an optimally preconditioned velocity-velocity matrix, it holds that:
\begin{equation}
  \begin{bmatrix} \mathbf{I} & \overline{\mathbf{A}}_{vp}\\ \overline{\mathbf{A}}_{vp}^T & \mathbf{0} \end{bmatrix} \left( \begin{array}{c} \mathbf{y}_{u,i} \\ \mathbf{y}_{p,i}  \end{array}  \right) = \lambda_i  \left( \begin{array}{c} \mathbf{y}_{u,i} \\ \mathbf{y}_{p,i}  \end{array}  \right),
\end{equation}
where $\overline{\mathbf{A}}_{vp}=\mathbf{S}_u^\frac12 \mathbf{A}_{vp} \mathbf{S}_p^{\frac12,T}$ is the preconditioned velocity-pressure matrix. By block-wise decomposition of this eigenvalue problem, the eigenvalues,  $\lambda_i$, are obtained as:
\begin{equation}
1, \quad \frac{1 + \sqrt{1+4 \mu_i}}{2}, \quad \frac{1 - \sqrt{1+4 \mu_i}}{2},
\label{eq:eigenvalues}\end{equation}
with $\mu_i$ the non-negative eigenvalues of the positive (semi)definite matrix:
\begin{equation}
\overline{\mathbf{A}}_{vp}^{T} \overline{\mathbf{A}}_{vp}.
\label{eq:pressureMatrix}\end{equation}
The multiplicity of the eigenvalue $1$ is equal to the number of velocity functions minus the rank of the matrix in \eqref{eq:pressureMatrix}. Assuming this matrix is nonsingular (\emph{i.e.,}\ the pair of function spaces is inf-sup stable) this is equal to the number of velocity functions minus the number of pressure functions. The multiplicities of the eigenvalues $(1 + \sqrt{1+4 \mu_i})/2$ and $(1 - \sqrt{1+4 \mu_i})/2$ are equal to the multiplicity of the eigenvalue $\mu_i$.

In view of the symmetry of $\overline{\mathbf{A}}$, the condition number of the preconditioned matrix $\overline{\mathbf{A}}$ corresponds to the ratio between the (in magnitude) largest and smallest eigenvalue. Equation~\eqref{eq:eigenvalues} relates the eigenvalue spectrum of this matrix to that of the positive (semi)definite matrix \eqref{eq:pressureMatrix}. The spectrum of the matrix \eqref{eq:pressureMatrix} depends on the preconditioning matrices $\mathbf{S}_u^\frac12$ and $\mathbf{S}_p^\frac12$, the former of which has already been fixed to optimally precondition the velocity-velocity matrix. The preconditioning matrix $\mathbf{S}_p^\frac12$ is optimal if the resulting condition number of $\overline{\mathbf{A}}$ is minimal. To derive the optimal preconditioner $\mathbf{S}_p^\frac12$ we first note that:
\begin{equation}
 \frac{|\lambda|_{\rm max} \left( \mu_1, \mu_2, ... , \mu_{n_p} \right)}{|\lambda|_{\rm min} \left( \mu_1, \mu_2, ... , \mu_{n_p} \right)} \geq \frac{|\lambda|_{\rm max} \left( \mu_1 \right)}{|\lambda|_{\rm min} \left( \mu_1 \right)},
\end{equation}
such that optimality requires $\mu_1 = \mu_2 = ... = \mu_{n_p} \equiv \mu$. Formal optimization of $\mu$ yields $\mu=2$ resulting in:
\begin{equation}
\kappa_2(\overline{\mathbf{A}})= \frac{ \left| \frac{1 + \sqrt{1+4 \mu}}{2} \right| }{ \left| \frac{1 - \sqrt{1+4 \mu}}{2} \right| } = \frac{ \left| 2 \right| }{ \left| -1 \right| } = 2.
\label{eq:optimalmmixed}
\end{equation}
All eigenvalues of $\overline{\mathbf{A}}_{vp}^T \overline{\mathbf{A}}_{vp}$ are equal to $\mu=2$ if the matrix  is equal to $2\mathbf{I}$, such that:
\begin{equation}
 \overline{\mathbf{A}}_{vp}^{T} \overline{\mathbf{A}}_{vp} = \left( \mathbf{S}_u^\frac12 \mathbf{A}_{vp} \mathbf{S}_p^{\frac12,T} \right)^T \left( \mathbf{S}_u^\frac12 \mathbf{A}_{vp} \mathbf{S}_p^{\frac12,T} \right) = \mathbf{S}_p^\frac12 \left( \mathbf{A}_{vp}^{T} \mathbf{S}_{u} \mathbf{A}_{vp} \right) \mathbf{S}_p^{\frac12,T} = 2\mathbf{I}.
\end{equation}
From the last equality in this expression one then obtains an optimal preconditioning matrix, $\mathbf{S}_p^\frac12$, for the pressure block:
\begin{equation}
\mathbf{S}_p = \mathbf{S}_p^{\frac12,T}\mathbf{S}_p^\frac12 = \left( \tfrac{1}{2}\mathbf{A}_{vp}^{T} \mathbf{S}_{u} \mathbf{A}_{vp} \right)^{-1}.
\label{eq:schurr}\end{equation}

The procedure for determining the pressure space preconditioning matrix, $\mathbf{S}_p$, according to \eqref{eq:schurr}, is as follows.  One first computes $\mathbf{S}_u$ as the CbAS preconditioner of the velocity-velocity matrix $\mathbf{A}_{vu}$ for each spatial direction. Using this velocity preconditioner one then computes $\mathbf{S}_p$ as the CbAS preconditioner of the matrix $\frac12 \mathbf{A}_{vp}^{T} \mathbf{S}_{u} \mathbf{A}_{vp}$.

For the Navier-Stokes equations, $\mathbf{A}_{vu}$ and consequently $\mathbf{S}_u$ are not symmetric and for some stabilized variational forms also $\mathbf{A}_{qu} \neq \mathbf{A}_{vp}^T$. As mentioned in Section~\ref{sec:preconSingle}, one can still apply the Additive-Schwarz framework to the nonsymmetric matrix $\frac12\mathbf{A}_{qu} \mathbf{S}_{u} \mathbf{A}_{vp}$ for such problems.

We note that the use of separate preconditioning blocks for the Stokes problem (or other saddle point problems) has been explored before. It is a well-established concept to use an approximate inverse of $\mathbf{A}_{vu}$ for the velocity and to use an approximate inverse of the Schur complement for the pressure (see \emph{e.g.,}\ \cite{Murphy2000}). In practice the factor $2$ in \eqref{eq:optimalmmixed}-\eqref{eq:schurr} does not significantly influence the condition number or the convergence behavior of iterative methods for the preconditioned system and could be omitted to simplify the procedure. For the examples in Section~\ref{sec:results} this factor is retained, however.

\begin{remark}\label{remark:scaling}
Because the CbAS preconditioner is composed of a summation of local inverses, the preconditioner deviates from the inverse at entries corresponding to functions that appear in multiple blocks. More specifically, it holds that $\mathbf{S}_u\mathbf{A}_{vu} \approx {\rm diag}( k_{u,1}, k_{u,2}, \ldots, k_{u,n_u} )$, where $k_{u,j}\geq 1$ represents the number of blocks that contain the (velocity) basis function with index $j$. The preconditioner can hence be further improved by diagonal scaling according to:
\begin{equation}
\mathbf{S}_{u}^{\rm scaled} =  {\rm diag}( k_{u,1}^{-\frac{1}{2}}, k_{u,2}^{-\frac{1}{2}}, \ldots, k_{u,n_u}^{-\frac{1}{2}} ) \, \mathbf{S}_{u} \,  {\rm diag}( k_{u,1}^{-\frac{1}{2}}, k_{u,2}^{-\frac{1}{2}}, \ldots, k_{u,n_u}^{-\frac{1}{2}} )
\label{eq:rescaling1}
\end{equation}
and, similarly for the pressure functions:
\begin{equation}
\mathbf{S}_{p}^{\rm scaled} =  {\rm diag}( k_{p,1}^{-\frac{1}{2}}, k_{p,2}^{-\frac{1}{2}}, \ldots, k_{p,n_p}^{-\frac{1}{2}} ) \, \mathbf{S}_{p} \,  {\rm diag}( k_{p,1}^{-\frac{1}{2}}, k_{p,2}^{-\frac{1}{2}}, \ldots, k_{p,n_p}^{-\frac{1}{2}} ),
\label{eq:rescaling2}
\end{equation}
with $\mathbf{S}_{p}$ the CbAS preconditioner of $\frac12\mathbf{A}_{qu}\mathbf{S}_{u}^{\rm scaled}\mathbf{A}_{vp}$. In the examples in Section~\ref{sec:results} we did, however, not observe a significant effect of rescaling \eqref{eq:rescaling1} and \eqref{eq:rescaling2} on the condition number.
\end{remark}\label{sec:preconMixed}

\section{Numerical examples}
In this section we verify the effectiveness of the CbAS preconditioner proposed in Section~\ref{sec:preconIntro}. We consider four problems to test various aspects of the preconditioner: a Poisson problem with boundary conditions imposed by the symmetric and nonsymmetric Nitsche method in Section~\ref{sec:Poisson}; an SUPG-stabilized convection-dominated convection-diffusion problem with the symmetric Nitsche method in Section~\ref{sec:SUPG}; a Stokes flow problem with the symmetric Nitsche method in Section~\ref{sec:Stokes}; and the steady and transient incompressible Navier-Stokes equations with the symmetric Nitsche method in Section~\ref{sec:NavierStokes}. These problems range from a symmetric positive definite (SPD) single-field problem for which SIPIC is suitable, to nonsymmetric indefinite mixed problems which are well outside the scope of SIPIC.

\begin{figure}
  \begin{center}
  \includegraphics[width=0.5\textwidth]{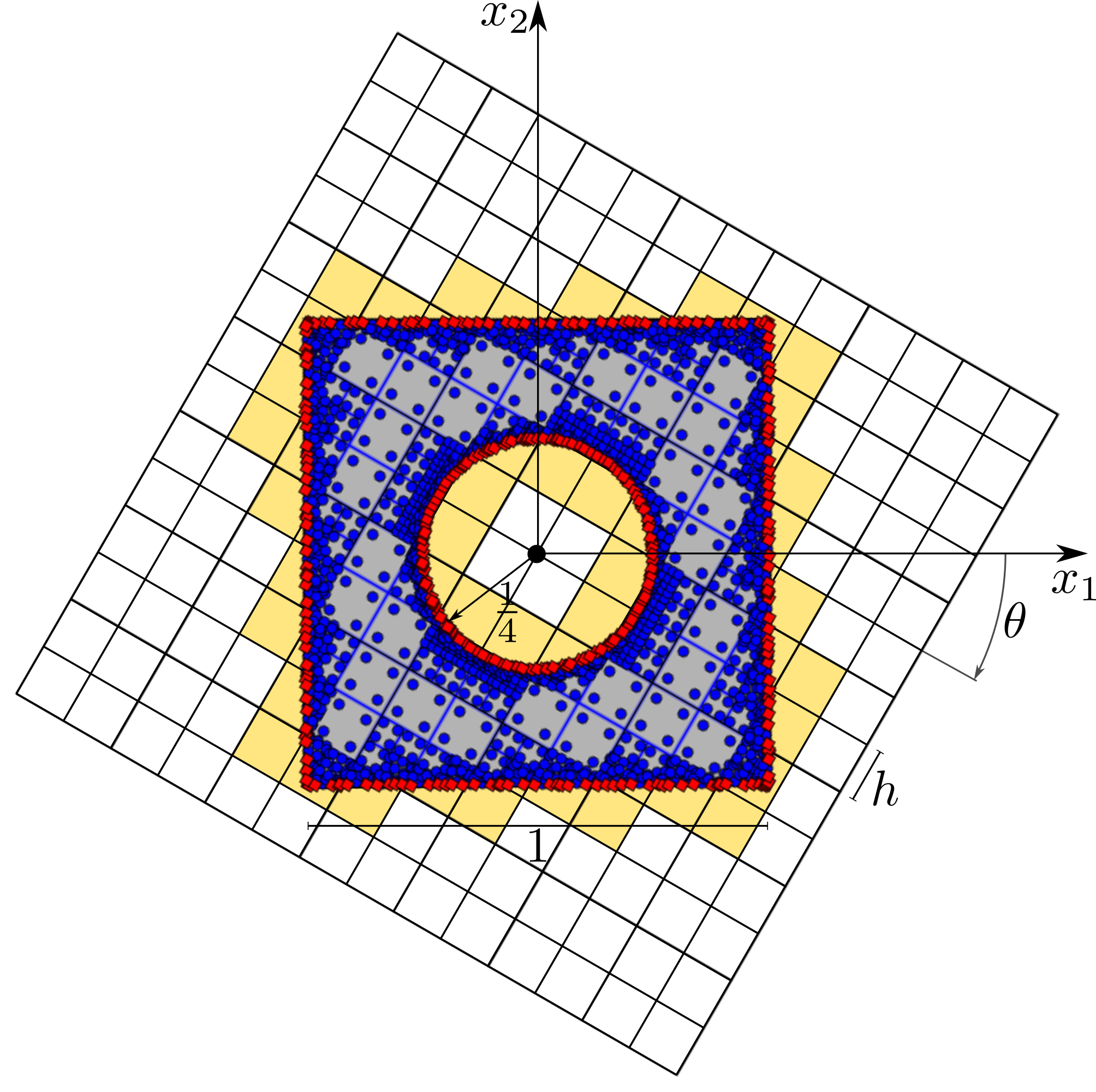}
  \caption{Schematic representation of the physical domain $\Omega$ embedded in a simple rectilinear encapsulating domain. The trimmed elements that intersect the physical domain are indicated in yellow and are used in the computation. The white elements do not intersect the domain and are not considered. The blue circles indicate volumetric integration points and the red squares indicate the Gauss points for boundary integration. Note that for graphical clarity the grid size in this figure is not to scale.\label{fig:GaussPoints}}
  \end{center}
\end{figure}

Except for the transient Navier-Stokes problem in Section~\ref{sec:NSvalidation}, all problems are posed on the same physical domain. This domain $\Omega$ consists of an origin-centered unit square $\left(-\frac12,\frac12\right)^2$ with, also at the origin, a circular exclusion of radius $\frac14$:
\begin{equation}
\Omega = \bigg\{ x = \left(x_1,x_2\right) : |x|_\infty < \frac12 \hspace{.1cm} , \hspace{.1cm} |x|_e > \frac14 \bigg\}. 
\end{equation}
By convention we denote the horizontal direction by $x_1$ and the vertical direction by $x_2$. The encapsulating domain is discretized by a uniform tensor product mesh with mesh size $h=\frac1{16}$. The grid has a vertex at the origin and is initially aligned with the edges of the physical domain (\emph{i.e.,}\ $\theta = \SI{0}{\degree}$). The encapsulating domain (and consequently the grid) is then rotated in $100$ steps until $\theta = \SI{45}{\degree}$. This generates different discretizations of the same domain with the same grid size but with varying smallest volume fractions $\eta$. The physical domain, the encapsulating domain and the grid are depicted in Figure~\ref{fig:GaussPoints}. For all different arrangements we assemble the full system matrix and compute the original and preconditioned condition number $\kappa_2$ or eigenvalue ratio $\rho$ (see Section~\ref{sec:conditioning}). These values provide good indications of the complexity of the system for either the Conjugate Gradient or the Generalized Minimal RESidual iterative solver.

We consider discrete approximations based on isogeometric $C^1(\Omega)$-continuous bivariate B-spline bases, $\mathcal{S}^2_1(\Omega)$, for the scalar single-field problems. For the two-field mixed problems, we consider the Taylor-Hood pair comprising the velocity space $\mathcal{S}^2_0(\Omega) \times \mathcal{S}^2_0(\Omega)$ and the pressure space $\mathcal{S}^1_0(\Omega)$. The geometry is approximated by the bisection-based tessellation scheme proposed in \cite{Verhoosel2015} with a maximal refinement depth of three. The number of Gauss points is taken such that all functionals are integrated exactly over the tessellated domain. For second order terms these points are displayed in Figure~\ref{fig:GaussPoints}. Note that variational forms with second order bases as used here contain fourth order terms and require more integration points than indicated in the figure. The condition numbers and eigenvalues are computed by a power algorithm that is terminated when the relative difference between subsequent iterations is smaller than $10^{-6}$.
\label{sec:results}
\subsection{Poisson's problem with the symmetric and nonsymmetric Nitsche method}
We consider the Poisson problem:
\begin{equation}
 \begin{cases}
  -\Delta u = f = 1 & {\rm in \ } \Omega, \\
  u = g^D = 0 & {\rm on \ } \Gamma^D = \{|x|_\infty = \frac12\}, \\
  \partial_n u = g^N = 0 & {\rm on \ } \Gamma^N = \{|x|_e = \frac14\},
 \end{cases}
 \label{eq:Poissonstrong}
\end{equation}
with the symmetric \cite{Nitsche1971} and nonsymmetric \cite{Freund1995,Burman2012,Schillinger2016} variants of Nitsche's method to impose boundary conditions. The homogeneous Dirichlet condition is imposed on the boundary of the unit square and the homogeneous Neumann condition is imposed on the boundary of the circular exclusion.

The variational form of the problem with the symmetric Nitsche method, resulting in an SPD matrix, is:
\begin{equation}\left\{ \begin{array}{l}
\mbox{Find } u_h \in \mathcal{V}_h(\Omega) \mbox{ such that } \forall v_h \in \mathcal{V}_h(\Omega)\mbox{:} \\[0.1em]
\int_\Omega \nabla v_h \cdot \nabla u_h {\rm d}V + \int_{\Gamma^D} - u_h \partial_n v_h - v_h \partial_n u_h + \beta v_h u_h {\rm d}S \\[0.1em]
\quad = \int_\Omega v_h f {\rm d}V + \int_{\Gamma^D} - g^D \partial_n v_h  + \beta v_h g^D {\rm d}S + \int_{\Gamma^N} v_h g^N {\rm d}S.
\end{array}\right.
\label{eq:PoissonSymmetric}
\end{equation}
In this formulation $\beta$ denotes the element-wise stabilization parameter that is required to ensure coercivity of the the bilinear operator. A lower bound for $\beta$ on each element, $i$, is given by:
\begin{equation}
 \beta_i > C_i = \max_{v_h \in \mathcal{V}_h(\Omega)} \frac{\| n \cdot \nabla v_h \|_{L^2(\Gamma_i^D)}^2}{\| \nabla v_h \|_{L^2(\Omega_i^{\rm tr})}^2}.
\end{equation}
To compute $C_i$ we solve a local generalized eigenvalue problem following the approach in \cite{Embar2010} and we set $\beta_i = 2C_i$. Let us mention that to improve the conditioning of this generalized eigenvalue problem we perform a local change of basis, as described in \cite[Appendix A]{SIPIC}. For shape regular trimmed elements it holds that $C_i \sim |\Gamma_i^D|/|\Omega_i^{\rm tr}| \sim 1/\hat{h}_i$ with $\hat{h}_i$ a typical length scale associated with the trimmed element $\Omega_i^{\rm tr}$. The variational form of the problem with the nonsymmetric Nitsche method, resulting in a nonsymmetric positive definite matrix, is:
\begin{equation}\left\{ \begin{array}{l}
\mbox{Find } u_h \in \mathcal{V}_h(\Omega) \mbox{ such that } \forall v_h \in \mathcal{V}_h(\Omega)\mbox{:} \\[0.1em]
\int_\Omega \nabla v_h \cdot \nabla u_h {\rm d}V + \int_{\Gamma^D} u_h \partial_n v_h - v_h \partial_n u_h + \frac{1}{h} v_h u_h {\rm d}S \\[0.1em]
\quad = \int_\Omega v_h f {\rm d}V + \int_{\Gamma^D} g^D \partial_n v_h + \frac{1}{h} v_h g^D {\rm d}S + \int_{\Gamma^N} v_h g^N {\rm d}S.
\end{array}\right.
\label{eq:PoissonNonsymmetric}
\end{equation}
Note that this nonsymmetric Nitsche method does not strictly require a stabilization parameter because the edge terms through $\Gamma^D$ cancel when $v_h=u_h$ \cite{Oden1998,Burman2012}. Both variational forms \eqref{eq:PoissonSymmetric} and \eqref{eq:PoissonNonsymmetric} are consistent with \eqref{eq:Poissonstrong}. On the full space $H^1(\Omega)$ both bilinear forms are not bounded and the bilinear form in \eqref{eq:PoissonSymmetric} is not coercive, see \emph{e.g.,}\ \cite{ErnGuermond}. Therefore \eqref{eq:PoissonSymmetric} and \eqref{eq:PoissonNonsymmetric} are posed on the finite dimensional space $\mathcal{V}_h(\Omega) \subset H^1(\Omega)$, on which boundedness and coercivity with respect to the $H^1$-norm are satisfied for both bilinear forms.

\begin{figure}
 \centering
 \begin{subfigure}[t]{.49\textwidth}
  \includegraphics[width=\textwidth]{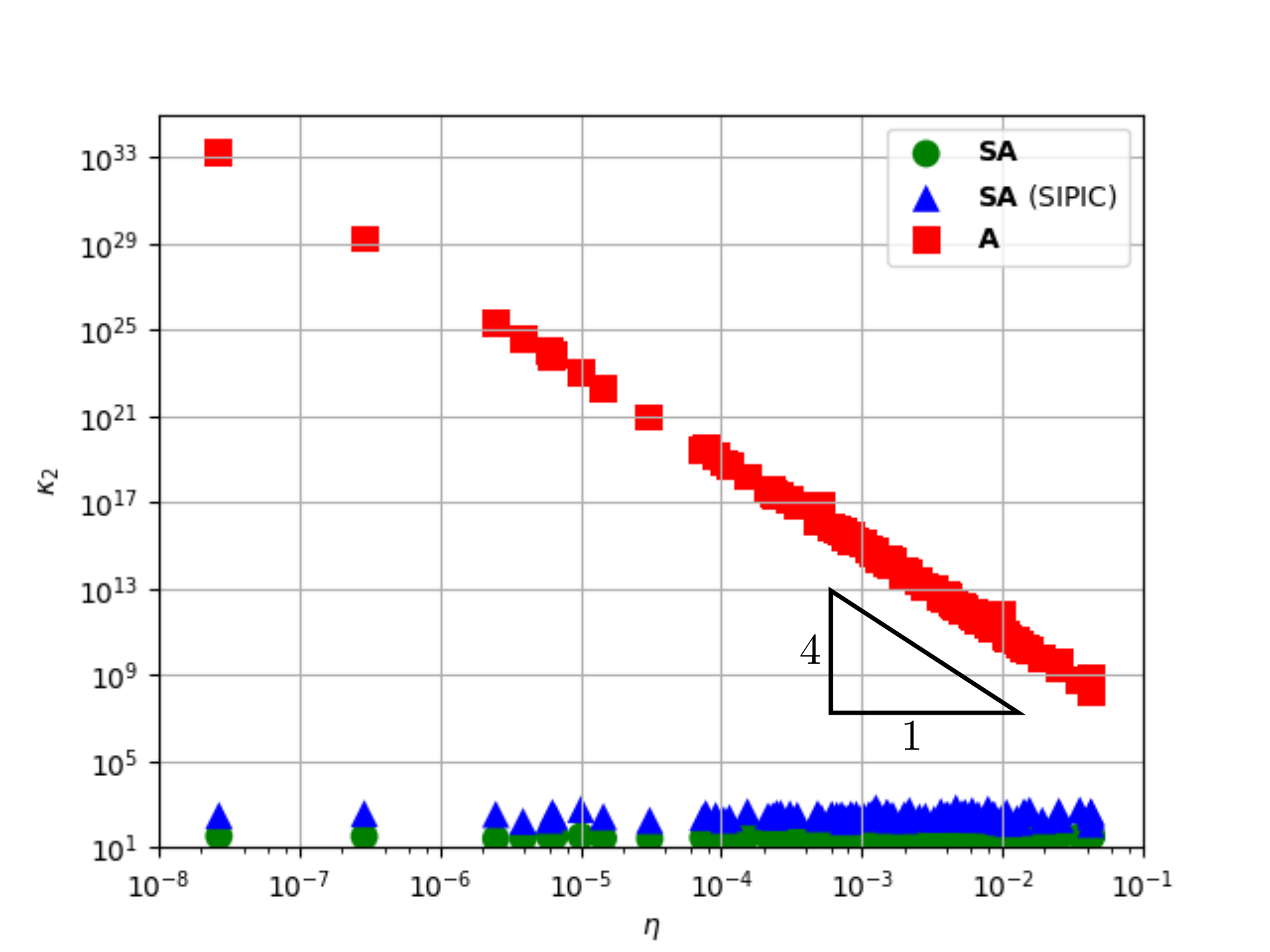}
  \caption{}
  \label{fig:PS_eta}
 \end{subfigure}
 \hfill
 \begin{subfigure}[t]{.49\textwidth}
  \includegraphics[width=\textwidth]{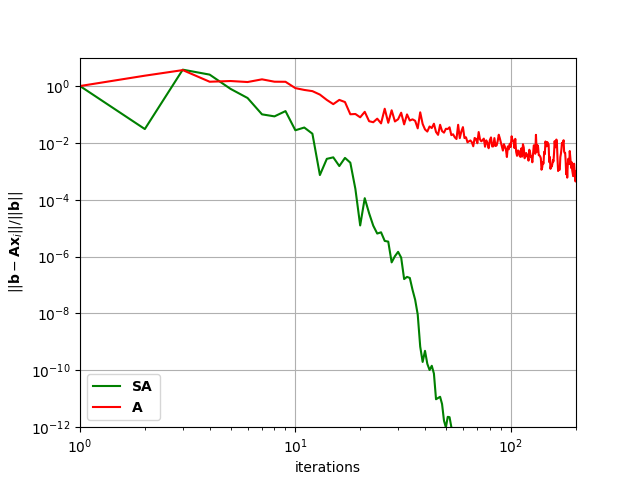}
  \caption{}
  \label{fig:PS_residual}
 \end{subfigure}
 \caption{(a) Original and preconditioned condition numbers (using SIPIC and CbAS) for the Poisson problem with the traditional symmetric Nitsche method as in \eqref{eq:PoissonSymmetric}. (b) The relative (preconditioned) residual error versus the number of CG iterations for an arrangement with $\theta = \SI{25}{\degree}$ yielding $\eta=9\cdot10^{-4}$. \label{fig:PS}}
\end{figure}

In Figure~\ref{fig:PS_eta} we present the condition numbers corresponding to the symmetric variational form \eqref{eq:PoissonSymmetric} with and without the CbAS preconditioner. We plot the condition numbers versus the smallest volume fraction for all arrangements. Note that the same test case was considered in \cite[Section 4.3]{SIPIC}. Because of the smooth second order B-spline basis and the SPD character of the matrices, the SIPIC preconditioner yields results similar to those of CbAS, which illustrates the interpretation of CbAS as a generalization of SIPIC. We observe that without preconditioning the condition number scales with the smallest volume fraction to the power $2p+1-2/d=4$ (see Section~\ref{sec:conditioning}). We also observe that with CbAS preconditioning the system is well-conditioned and robust with respect to the smallest volume fraction, in the sense that for all volume fractions we obtain condition numbers in the range of $24 \leq \kappa_2(\mathbf{A}) \leq 38$ and that a scaling relation is not observed. A typical CG solver convergence plot is shown in Figure~\ref{fig:PS_residual}, from which it is observed that CbAS preconditioning results in substantially improved convergence behavior of the iterative solver. 

\begin{figure}
 \centering
 \begin{subfigure}[t]{.49\textwidth}
  \includegraphics[width=\textwidth]{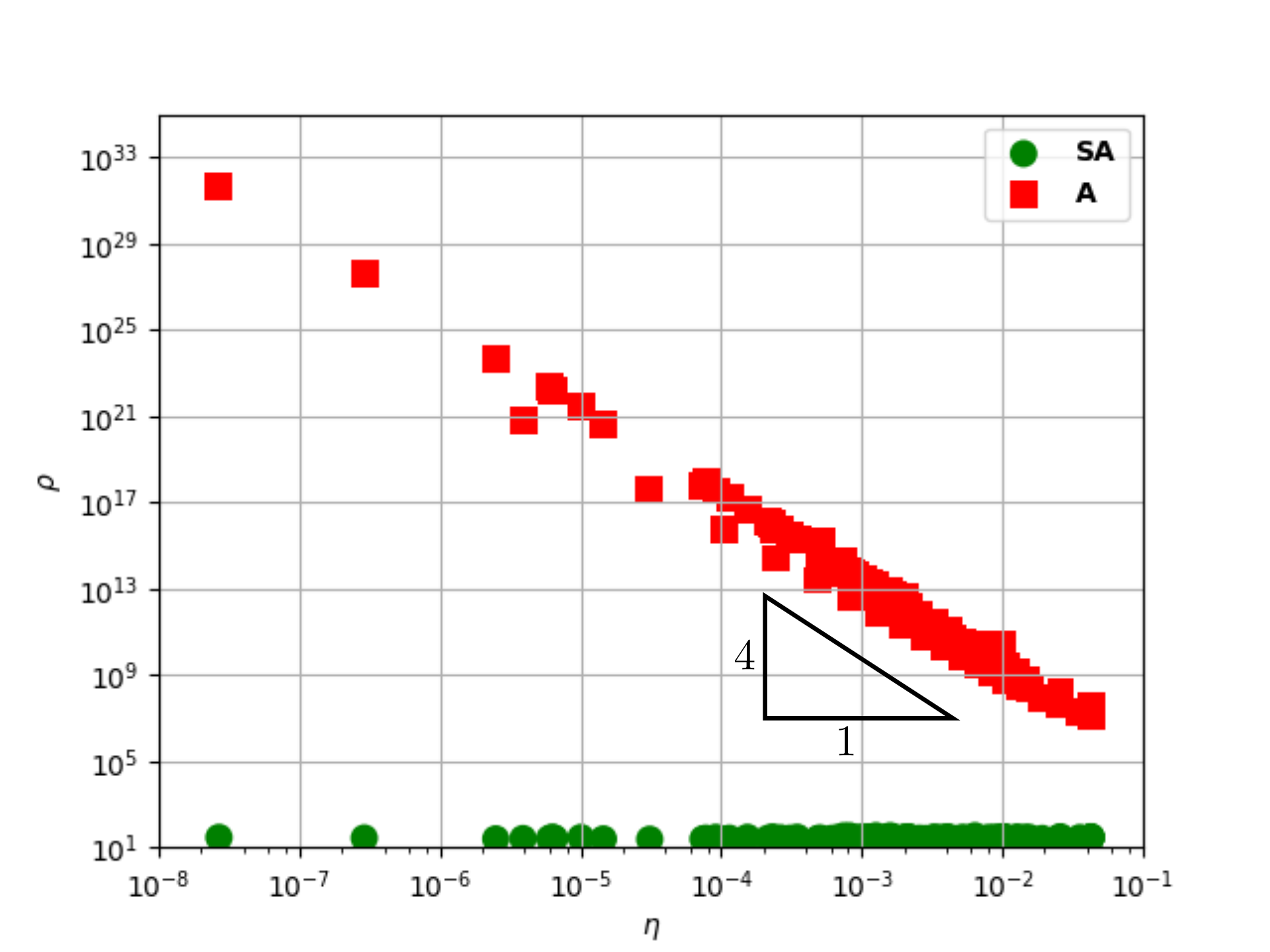}
  \caption{}
  \label{fig:PN_eta}
 \end{subfigure}
 \hfill
 \begin{subfigure}[t]{.49\textwidth}
  \includegraphics[width=\textwidth]{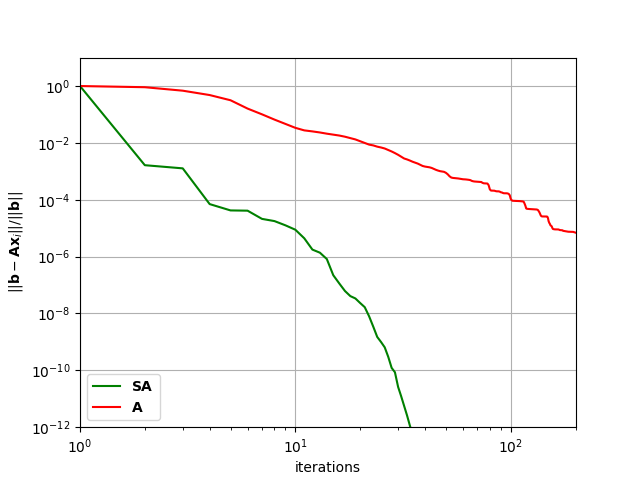}
  \caption{}
  \label{fig:PN_residual}
 \end{subfigure}
 \caption{(a) Original and preconditioned eigenvalue ratios for the Poisson problem with the nonsymmetric Nitsche method. (b) The relative (preconditioned) residual error versus the number of GMRES iterations for an arrangement with $\theta = \SI{25}{\degree}$ yielding $\eta=9\cdot10^{-4}$.\label{fig:PN}}
\end{figure}

The eigenvalue ratios with and without preconditioning of the nonsymmetric variational form in \eqref{eq:PoissonNonsymmetric} are plotted in Figure~\ref{fig:PN_eta}. This figure conveys that the CbAS preconditioner is also effective for this non-SPD system, as for all volume fractions we obtain eigenvalue ratios in the range of $23 \leq \rho(\mathbf{A}) \leq 34$. The results show that the conditioning with the nonsymmetric Nitsche method is almost identical to the conditioning with the symmetric Nitsche method, both with and without preconditioning. This can be explained by the detailed analysis of the ill-conditioning mechanism presented in \cite[Section~3.2]{SIPIC}. Here it is shown that on small cut elements, the contribution of the stabilization terms in the symmetric form are of the same magnitude as the other terms. Therefore we expect approximately the same norm of $\|\mathbf{A}^{-1}\|$ and consequently the same conditioning of the symmetric and the nonsymmetric forms. The effect of CbAS preconditioning on the convergence of the residual in a GMRES solver is illustrated in Figure~\ref{fig:PN_residual}.

\label{sec:Poisson}
\subsection{Convection-diffusion prolem}
\begin{figure}
  \begin{center}
  \includegraphics[width=0.65\textwidth]{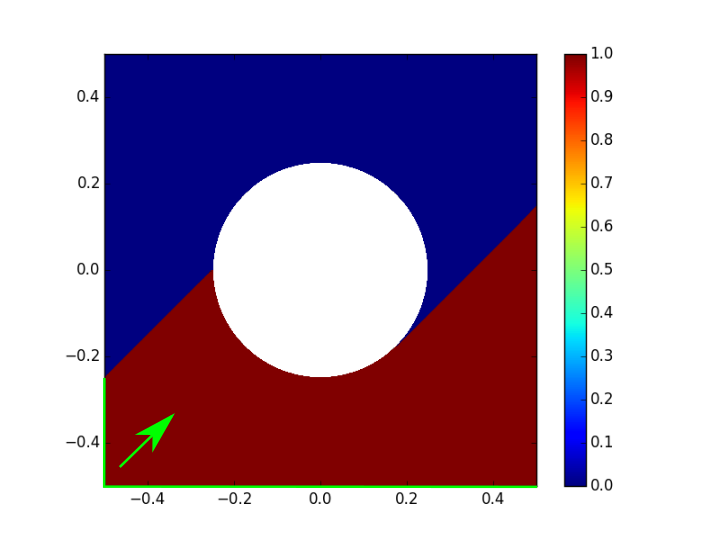}
  \caption{Limit of the solution to the convection-diffusion problem \eqref{eq:SUPGstrong} for $\varepsilon \to +0$. The non-homogeneous part of the Dirichlet boundary, $\Gamma_1^D$, is indicated in green, and the green arrow indicates the direction of the convective velocity $w$.\label{fig:SUPGSolution}}
  \end{center}
\end{figure}

We consider the convection-dominated convection-diffusion problem:
\begin{equation}
 \begin{cases}
  {\rm div}\left( wu - \varepsilon \nabla u \right) = 0 & {\rm in \ } \Omega, \\
  u = g^D & {\rm on \ } \Gamma^D = \partial\Omega,
 \end{cases}
 \label{eq:SUPGstrong}
\end{equation}
with field variable $u$, convective velocity $w = (1,1)$ and diffusion coefficient $0 < \varepsilon = 10^{-6} \ll 1$, such that $|w|_e \gg \varepsilon$ and the convection is clearly dominant. Because both the boundary of the unit square ($|x|_\infty = \frac12$) and the boundary of the circular exclusion ($|x|_e = \frac14$) are partially inflow and outflow boundaries, we pose Dirichlet conditions on the full boundary $\partial\Omega$. The function $g^D$ prescribes $u=1$ on the lower boundary ($x_2 = -\frac12$) and part of the left boundary ($x_1 = -\frac12$ and $x_2<-\frac14$) and prescribes $u=0$ on the remainder:
\begin{equation}
 g^D = 
 \begin{cases}
 1 \quad & \mbox{on } \Gamma_1^D = \big\{ \big(x_1,-\frac12\big) \big| x_1 \in \big(-\frac12,\frac12\big) \big\} \cup \big\{ \big(-\frac12,x_2\big) \big| x_2 \in \big(-\frac12,-\frac14\big) \big\}, \\
 0 \quad & \mbox{on } \Gamma^D \setminus \Gamma_1^D.
 \end{cases}
 \label{eq:supgdirichletdata}
\end{equation}
The function $g^D$ does not affect the system matrix and consequently does not affect the conditioning however. The solution to \eqref{eq:SUPGstrong}--\eqref{eq:supgdirichletdata} for the limit $\varepsilon \to +0$ is shown in Figure~\ref{fig:SUPGSolution}. Note that the boundary layers at the right boundary ($x_1 = \frac12$) and at the boundary of the circular exclusion ($|x|_e = \frac14$) of thickness $\varepsilon/|w|_e$ disappear for $\varepsilon \to +0$. For $\varepsilon = 0$ the boundary conditions are violated, such that a solution in $H^1(\Omega)$ no longer exists. This is evident as for $\varepsilon = 0$ the problem posed in \eqref{eq:SUPGstrong} is no longer elliptic. The solution in Figure~\ref{fig:SUPGSolution} is therefore only a limit in $L^2(\Omega)$, as for $\varepsilon \to +0$ the solutions diverge in $H^1(\Omega)$.

We employ the variational form introduced in \cite{Bazilevs2007} in which the convective terms are stabilized by Streamlined Upwind Petrov-Galerkin (SUPG) terms \cite{Brooks1982} and Dirichlet boundary conditions are imposed by Nitsche's method:
\begin{equation}\left\{ \begin{array}{l}
\mbox{Find } u_h \in \mathcal{V}_h(\Omega) \mbox{ such that } \forall v_h \in \mathcal{V}_h(\Omega)\mbox{:} \\[0.1em]
\int_\Omega \left( -u_h w \cdot \nabla v_h + \varepsilon \nabla v_h \cdot \nabla u_h + \tau w \cdot \nabla v_h {\rm div}\left( wu_h - \varepsilon \nabla u_h \right) \right) {\rm d}V \\[0.1em]
+ \int_{\partial\Omega} \left( {\rm max}\left(0,n \cdot w\right)v_h u_h - \varepsilon \left( v_h \partial_n u_h + u_h \partial_n v_h \right) + \varepsilon \beta v_h u_h \right) {\rm d}S \\[0.1em]
\quad = \int_{\partial\Omega} \left( - {\rm min}\left(0,n \cdot w\right)v_h g^D - \varepsilon g^D \partial_n v_h + \varepsilon \beta v_h g^D \right) {\rm d}S.
\end{array}\right.
\label{eq:SUPGweak}
\end{equation}
The same element-wise stabilization parameter $\beta$ as in \eqref{eq:PoissonSymmetric} is employed. Different choices for the SUPG parameter $\tau$ are motivated in \cite{Tezduyar2000}. In all our examples we use $\tau = h^*/\left(2|w|_e\right)$ with $h^*$ the maximal element length in the direction of velocity $w$. For uniform tensor product grids this implies $\tau = h/\left(2 \max_i \left(|w \cdot e_i|\right)\right)$, with $e_i$ the unit vector in the direction of a grid line. We therefore have:
\begin{equation}
\tau = \frac{h}{2 \sqrt{2} \sin\left( \frac14 \pi + \theta \right)}.
\label{eq:supgtau}
\end{equation}
Our computations based on $\tau$ according to \eqref{eq:supgtau} as a global parameter did not indicate a need to consider a local $\hat{h}_i^*$ on trimmed elements. We note that in the setting considered here all trimmed elements also experience the stabilizing effect of the Nitsche terms, and that this observation does not necessarily extend to other immersed implementations of SUPG.

\begin{figure}
 \centering
 \begin{subfigure}[t]{.49\textwidth}
  \includegraphics[width=\textwidth]{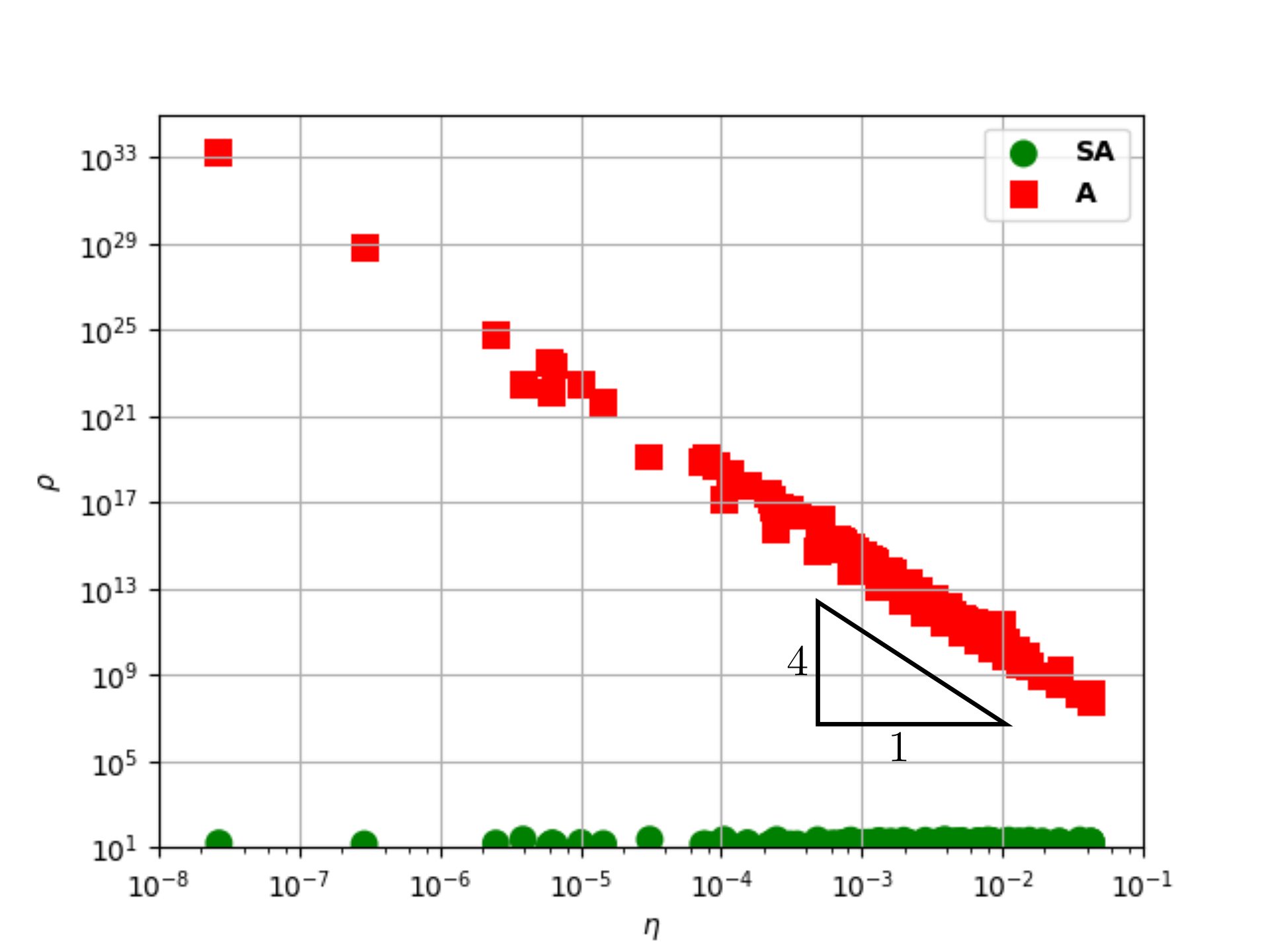}
  \caption{}
  \label{fig:SUPG_eta}
 \end{subfigure}
 \hfill
 \begin{subfigure}[t]{.49\textwidth}
  \includegraphics[width=\textwidth]{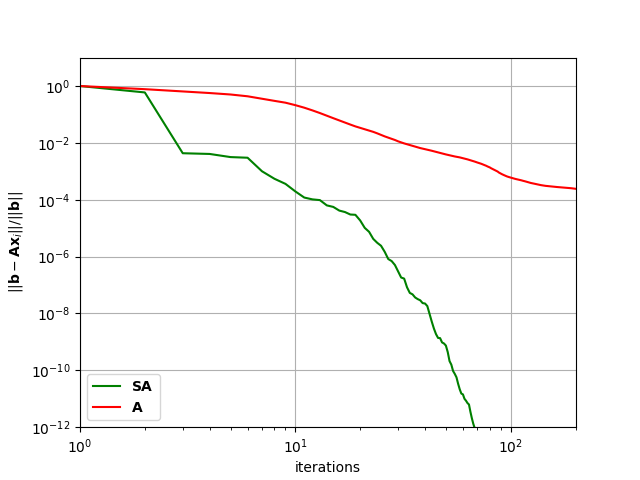}
  \caption{}
  \label{fig:SUPG_residual}
 \end{subfigure}
 \caption{(a) Original and preconditioned eigenvalue ratios for the convection-dominated convection-diffusion problem with boundary conditions imposed by Nitsche's method. (b) The relative (preconditioned) residual error versus the number of GMRES iterations for an arrangement with $\theta = \SI{25}{\degree}$ yielding $\eta=9\cdot10^{-4}$.
 \label{fig:SUPG}}
\end{figure}

The eigenvalue ratios of the resulting non-SPD matrices without preconditioning and with CbAS preconditioning are plotted in Figure~\ref{fig:SUPG_eta} for all arrangements. We again observe that without preconditioning the eigenvalue ratio scales with the smallest volume fraction to the power $4$. The results indicate that with preconditioning the system is well-conditioned and robust with respect to the smallest volume fraction, as all obtained eigenvalue ratios lie in the range $12 \leq \rho(\mathbf{A}) \leq 23$ and no scaling relation can be observed. The resulting positive effect on the GMRES solver convergence behavior is illustrated in Figure~\ref{fig:SUPG_residual}.

\label{sec:SUPG}
\subsection{Stokes problem}
\begin{figure}
 \centering
 \begin{subfigure}[t]{.49\textwidth}
  \includegraphics[width=\textwidth]{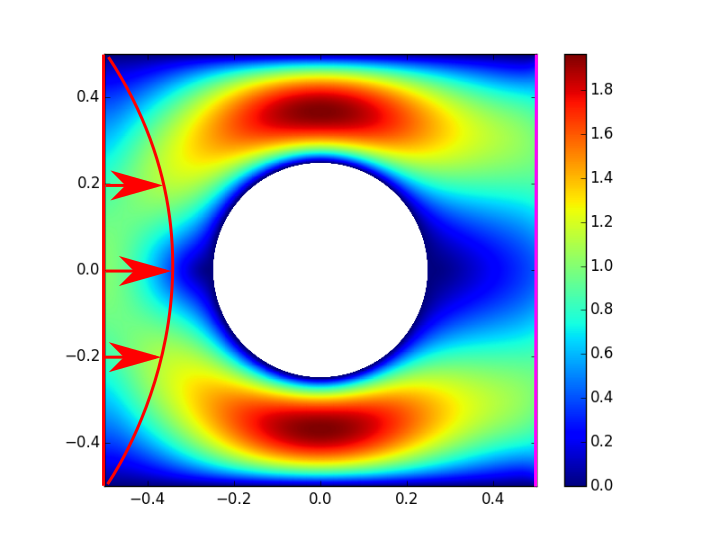}
  \caption{Velocity magnitude\label{fig:StokesSolutionVelocity}}
 \end{subfigure}
 \hfill
 \begin{subfigure}[t]{.49\textwidth}
  \includegraphics[width=\textwidth]{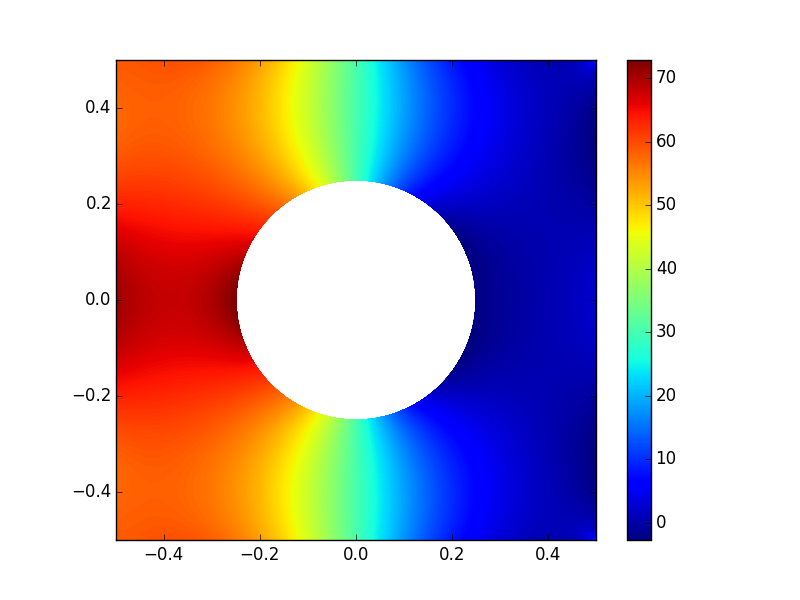}
  \caption{Pressure}
 \end{subfigure}
 \caption{Solution to the Stokes problem \eqref{eq:Stokesstrong} subject to the boundary conditions defined by \eqref{eq:gDStokes}. The inflow boundary $\Gamma_{\rm in}^D$ and the outflow boundary $\Gamma^N$ are indicated in Figure~\ref{fig:StokesSolutionVelocity} in respectively red and magenta.\label{fig:StokesSolution}}
\end{figure}

A Stokes flow problem is considered according to:
\begin{equation}
 \begin{cases}
  -{\rm div}\left(\nabla^s u - p \mathcal{I} \right) = 0 & {\rm in \ } \Omega, \\
  -{\rm div}\left(u\right) = 0 & {\rm in \ } \Omega, \\
  u = g^D & {\rm on \ } \Gamma^D = \partial\Omega \setminus \Gamma^N, \\
  n \cdot \nabla^s u - n p = g^N = 0 & {\rm on \ } \Gamma^N = \big\{ \big(\frac12,x_2\big) \big| x_2 \in \big( -\frac12,\frac12 \big) \big\}.     
 \end{cases}
 \label{eq:Stokesstrong}
\end{equation}
In this formulation $u$ and $p$ denote respectively the velocity and the pressure, $\nabla^s$ denotes the symmetric gradient operator and $\mathcal{I}$ denotes the second order identity tensor. We impose a Poiseuille profile on the left (inflow) boundary ($x_1=-\frac12$) and homogeneous Dirichlet conditions (no-slip) on the lower and upper boundary ($|x_2|=\frac12$) and on the boundary of the exclusion ($|x|_e=\frac14$):
\begin{equation}
 g^D = 
 \begin{cases}
 \left( 1-4x_2^2,0\right) \quad & \mbox{on } \Gamma_{\rm in}^D = \big\{ \big(-\frac12,x_2\big) \big| x_2 \in \big(-\frac12,\frac12\big) \big\}, \\
 \left(0,0\right) \quad & \mbox{on } \partial\Omega \setminus \left( \Gamma_{\rm in}^D \cup \Gamma^N \right),
 \end{cases}
\label{eq:gDStokes}\end{equation}
Moreover we impose homogeneous Neumann (traction free) conditions on the right (outflow) boundary $\Gamma^N$ ($x_1=\frac12$). The boundary condition data does not affect the conditioning. The solution to \eqref{eq:Stokesstrong} subject to these boundary conditions is plotted in Figure~\ref{fig:StokesSolution}.

The symmetric variational form with boundary conditions imposed by means of the symmetric Nitsche method is:
\begin{equation}\left\{ \begin{array}{l}
\mbox{Find } (u_h,p_h) \in \mathcal{V}_h(\Omega)\times \mathcal{Q}_h(\Omega) \mbox{ such that:} \\[0.1em]
\int_\Omega \nabla^s v_h : \nabla^s u_h {\rm d}V + \int_{\Gamma^D} - u_h \cdot n \cdot \nabla^s v_h - v_h \cdot n \cdot \nabla^s u_h + \beta v_h \cdot u_h {\rm d}S + \\[0.1em]
\int_\Omega - p_h {\rm div}(v_h) {\rm d}V + \int_{\Gamma^D} p_hn\cdot v_h {\rm d}S = \int_{\Gamma^D} - g^D \cdot n \cdot \nabla^s v_h + \beta v_h \cdot g^D {\rm d}S \hspace{.14cm} \forall v_h \in \mathcal{V}_h(\Omega), \\[0.1em]
\int_\Omega - q_h {\rm div}(u_h) {\rm d}V + \int_{\Gamma^D} q_hn\cdot u_h {\rm d}S = \int_{\Gamma^D} q_hn\cdot g^D {\rm d}S \hspace{2.78cm} \forall q_h \in \mathcal{Q}_h(\Omega).
\end{array}\right.
\label{eq:Stokesweak}
\end{equation}
The element-wise stabilization constant $\beta$ is again set to $\beta_i = 2C_i$, with $C_i$ here defined as:
\begin{equation}
 C_i = \max_{v_h \in \mathcal{V}_h(\Omega)} \frac{\| n \cdot \nabla^s v_h \|_{L^2(\Gamma_i^D)}^2}{\| \nabla^s v_h \|_{L^2(\Omega_i^{\rm tr})}^2}.
\label{eq:betaStokes}\end{equation}
A proper choice for the spaces $\mathcal{V}_h(\Omega)$ and $\mathcal{Q}_h(\Omega)$, such as the Taylor-Hood elements that we apply here, leads to inf-sup stability of \eqref{eq:Stokesweak}. An analysis of the stability and accuracy of different pairs of function spaces for this immersed problem is presented in \cite{Hoang2017}.

\begin{figure}
 \centering
 \begin{subfigure}[t]{.49\textwidth}
  \includegraphics[width=\textwidth]{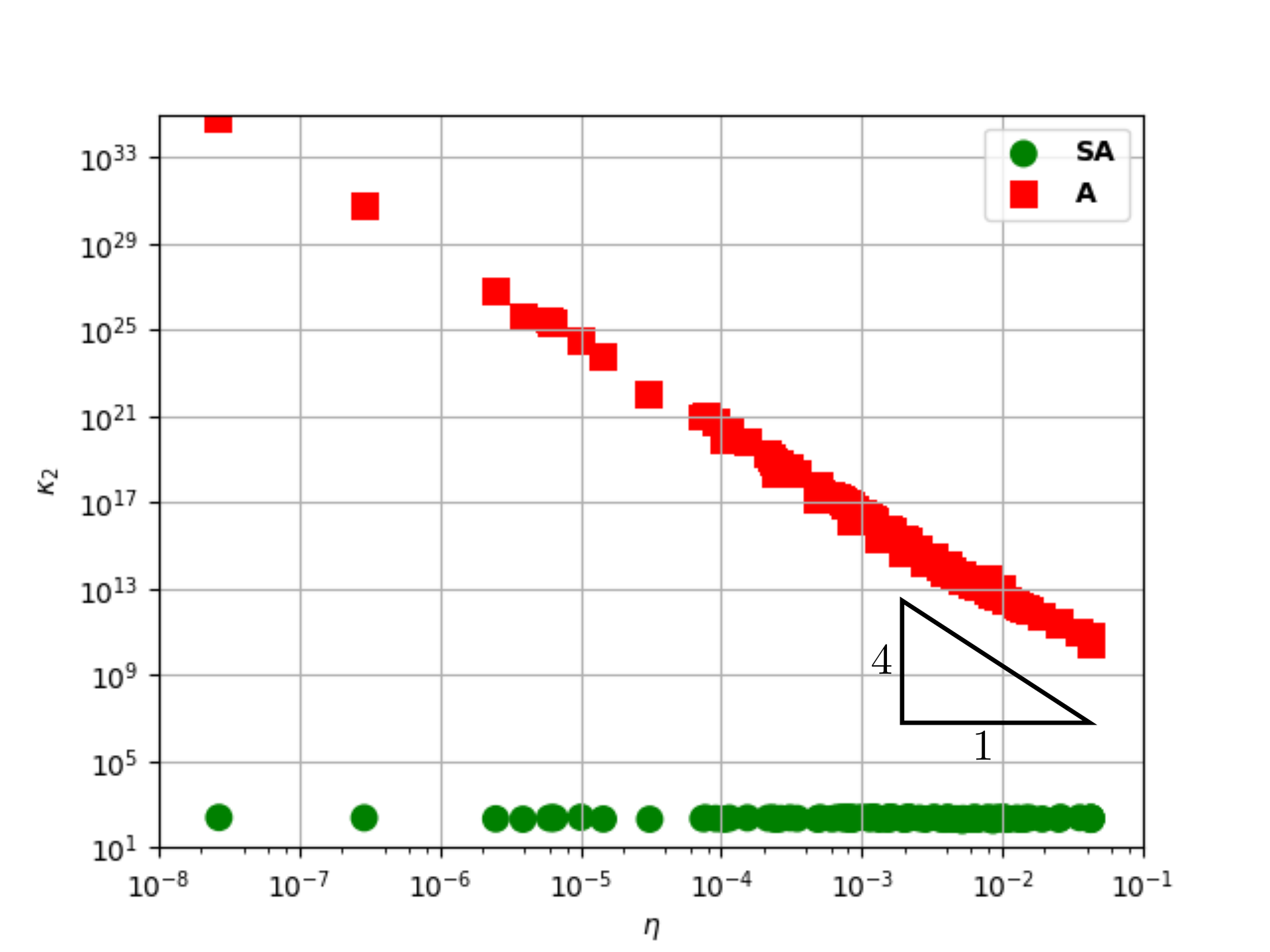}
  \caption{}
  \label{fig:Stokes_eta}
 \end{subfigure}
 \hfill
 \begin{subfigure}[t]{.49\textwidth}
  \includegraphics[width=\textwidth]{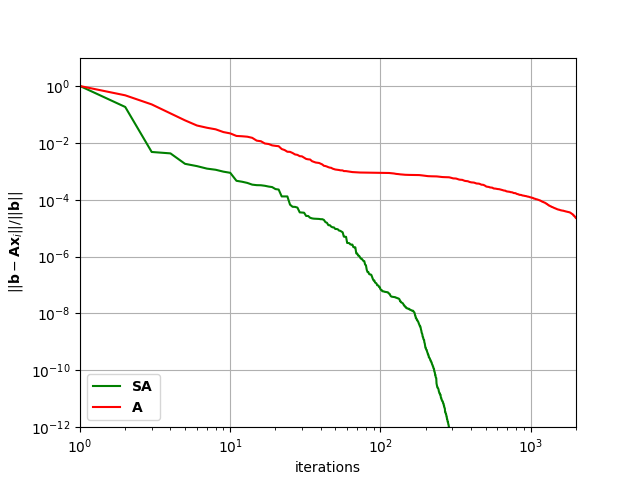}
  \caption{}
  \label{fig:Stokes_residual}
 \end{subfigure}
 \caption{(a) Original and preconditioned condition numbers for the Stokes problem with boundary conditions imposed by Nitsche's method. (b) The relative (preconditioned) residual error versus the number of GMRES iterations for an arrangement with $\theta = \SI{25}{\degree}$ yielding $\eta=9\cdot10^{-4}$. \label{fig:Stokes}}
\end{figure}

The resulting symmetric but indefinite matrices are preconditioned by the CbAS preconditioner using the Schur complement as described in Section~\ref{sec:preconMixed}. The condition numbers with and without preconditioning are presented in Figure~\ref{fig:Stokes_eta}. We observe that, also for this two-field mixed problem, the condition number without preconditioning scales with the smallest volume fraction to the power $4$. The CbAS preconditioner again provides a well-conditioned system that is robust with respect to the smallest volume fraction, with condition numbers ranging between $176$ and $247$. Figure~\ref{fig:Stokes_residual} illustrates the improved GMRES convergence behavior as a result of CbAS preconditioning.

\label{sec:Stokes}
\subsection{Navier-Stokes problem}\label{sec:NavierStokes}
In this section we demonstrate the applicability of the CbAS preconditioner to immersed flow problems. We consider both a steady and a transient test case.

\subsubsection{Steady Navier-Stokes problem}
\begin{figure}
 \centering
 \begin{subfigure}[t]{.49\textwidth}
  \includegraphics[width=\textwidth]{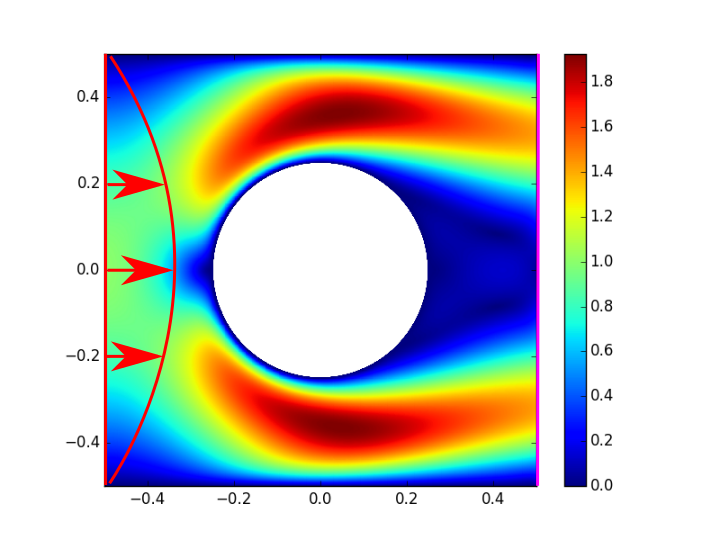}
  \caption{Velocity magnitude\label{fig:NavierStokesSolutionVelocity}}
 \end{subfigure}
 \hfill
 \begin{subfigure}[t]{.49\textwidth}
  \includegraphics[width=\textwidth]{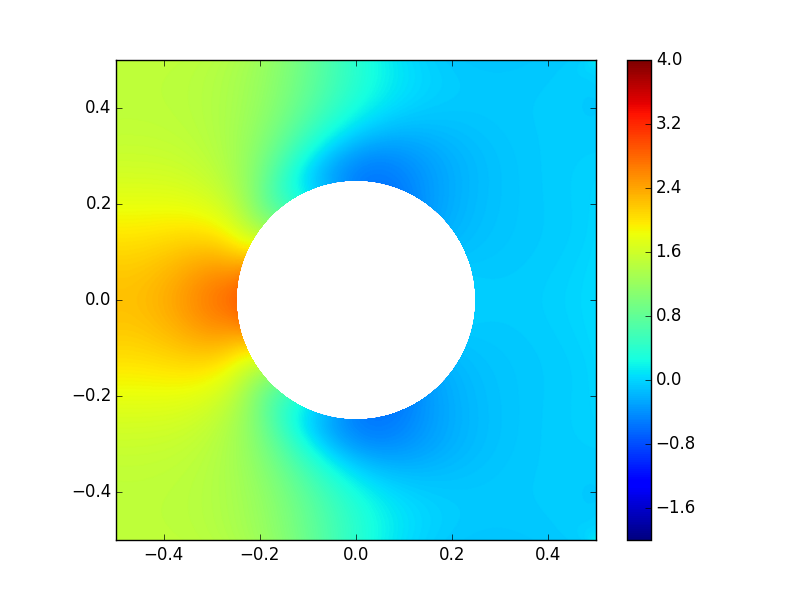}
  \caption{Pressure}
 \end{subfigure}
 \caption{Solution to the Navier-Stokes problem \eqref{eq:NSStrong} subject to the same boundary condition as the Stokes problem in Section~\ref{sec:Stokes}. The inflow boundary $\Gamma_{\rm in}^D$ and the outflow boundary $\Gamma^N$ are indicated in Figure~\ref{fig:NavierStokesSolutionVelocity} in respectively red and magenta.\label{fig:NSSolution}}
\end{figure}

To elaborate the considered Navier-Stokes problem, we first consider the following boundary value problem corresponding to the steady Oseen equations:
\begin{equation}
\begin{cases}
 {\rm div}\left( w \otimes u - 2 \nu \nabla^s u + p \mathcal{I} \right) = 0 \quad {\rm in\ } \Omega, \\
 -{\rm div}\left( u \right) = 0 \quad {\rm in\ } \Omega, \\
 u = g^D \quad {\rm on\ } \Gamma^D, \\
 2 \nu n \cdot \nabla^s u - n p  - \frac12 {\rm min}\left( 0, n \cdot w\right) u = g^N \quad {\rm on\ } \Gamma^N.
\end{cases}
\label{eq:NSStrong}\end{equation}
The dynamic viscosity is set to $\nu = 10^{-2} > 0$. Note that we add a directional do-nothing term to the Neumann condition to ensure well-posedness in case of backflow through $\Gamma^N$ \cite{Bruneau1996,Braack2014,Arndt2016}. By replacing $w$ in the convective term by $u$ in \eqref{eq:NSStrong} we recover the Navier-Stokes equations. We solve the nonlinear Navier-Stokes equations by means of a standard Picard iteration procedure, in which we solve the Oseen problem with $w^{n}=u^{n-1}$ to obtain $u^{n}$ in each iteration. As a stopping criterion we use:
\begin{equation}
\sqrt{\frac{\|u^{n+1} - u^n\|_{H^1(\Omega)}^2 + \|p^{n+1} - p^n\|_{L^2(\Omega)}^2}{\|u^{n+1} + u^n\|_{H^1(\Omega)}^2 + \|p^{n+1} + p^n\|_{L^2(\Omega)}^2}} \leq \mbox{tol} = 10^{-6},
\label{eq:stopcrit}\end{equation}
which is achieved in approximately 16 iterations for all arrangements. For the initial convective field $w^0$ we employ the solution to the Stokes problem. The solution to the steady Navier-Stokes problem is plotted in Figure~\ref{fig:NSSolution}.

To derive a suitable weak form of \eqref{eq:NSStrong}, we consider different parts of the Oseen equations separately. We consider a skew-symmetric form for the convection part \cite{Layton2008,Arndt2016}; based on the following identities:
\begin{equation}\begin{aligned}
& \int_\Omega v \cdot {\rm div}\left( w \otimes u \right) {\rm d}V \\
& = \frac12 \int_\Omega v \cdot {\rm div}\left( w \otimes u \right) {\rm d}V + \frac12 \int_\Omega {\rm div}\left( w \otimes u \cdot v \right) - u \otimes w : \nabla v {\rm d}V \\
& = \frac12 \int_\Omega \left(v \cdot u\right) \cancel{{\rm div}\left( w \right)} + v \otimes w : \nabla u - u \otimes w : \nabla v {\rm d}V + \frac12 \int_{\partial\Omega} \left( n \cdot w \right) \left( v \cdot u \right) {\rm d}S \\
& = \frac12 \int_\Omega v \otimes w : \nabla u - u \otimes w : \nabla v {\rm d}V + \frac12 \int_{\Gamma^D} \left( n \cdot \cancel{w} g^D \right) \left( v \cdot u \right) {\rm d}S \\
& \quad + \frac12 \int_{\Gamma^N} {\rm max}\left( 0, n \cdot w\right) \left( v \cdot u \right) {\rm d}S + \frac12 \int_{\Gamma^N} {\rm min}\left( 0, n \cdot w\right) \left( v \cdot u \right) {\rm d}S.
\end{aligned}\label{eq:NSconvection}\end{equation}
Note that the volume term in the ultimate expression is skew-symmetric in $u$ and $v$. We also mention that on the Dirichlet boundary the convective velocity $w$ is augmented with the prescribed velocity $g^D$ in order to reduce the nonlinearity in Navier-Stokes and improve the convergence speed of the Picard iteration. This is possible because our convective velocity field $w$ satisfies the Dirichlet conditions at all times. The integral over $\Gamma^D$ generally does not yield a positive semidefinite matrix because $n \cdot g^D < 0$ on inflow boundaries. This is repaired by the terms that impose the Dirichlet conditions, as will be shown in \eqref{eq:NSNitsche}. The first integral over $\Gamma^N$ yields a symmetric positive semidefinite matrix. The second integral over $\Gamma^N$ yields a negative semidefinite matrix, but this is canceled by the directional do-nothing part of the Neumann condition. The weak form of the stress contribution is:
\begin{equation}\begin{aligned}
& \int_\Omega v \cdot - {\rm div}\left( 2 \nu \nabla^s u - p \mathcal{I} \right) {\rm d}V = \\
& = \int_\Omega {\rm div}\left( - 2 \nu v \cdot \nabla^s u + pv \right) + 2 \nu \nabla^s v : \nabla^s u - p {\rm div}\left( v \right) {\rm d}V \\
& = \int_\Omega 2 \nu \nabla^s v : \nabla^s u - p {\rm div}\left( v \right) {\rm d}V + \int_{\partial\Omega} v \cdot - \left( 2 \nu n \cdot \nabla^s u - n p \right) {\rm d}S.
\end{aligned}\label{eq:NSstress}\end{equation}
The volume integral in this formulation is symmetrized by the weak incompressibility condition. The boundary integral in \eqref{eq:NSstress} is canceled on the Neumann boundary. The Dirichlet part of the boundary condition is made symmetric and coercive by the symmetric Nitsche method to impose Dirichlet conditions:
\begin{equation}\begin{aligned}
\int_{\Gamma^D} u & \cdot - \left( 2 \nu n \cdot \nabla^s v - nq \right) + \left( 2 \nu \beta - \frac12 {\rm min}\left( 0, n \cdot \cancel{w} g^D \right) \right) \left( v \cdot u \right) {\rm d}S \\
= \int_{\Gamma^D} g^D & \cdot - \left( 2 \nu n \cdot \nabla^s v - nq \right) + \left( 2 \nu \beta - \frac12 {\rm min}\left( 0, n \cdot \cancel{w} g^D \right) \right) \left( v \cdot g^D \right) {\rm d}S.
\end{aligned}\label{eq:NSNitsche}\end{equation}
We use the element-wise stabilization parameter $\beta$ as for the Stokes problem in Section~\ref{sec:Stokes}. Note that the convective term in \eqref{eq:NSNitsche} is only imposed on the inflow boundaries in order to stabilize the formulation, see \emph{e.g.,}\ \cite{Bazilevs2007}. As in \eqref{eq:NSconvection} the convective terms are again augmented.

When we assemble the terms above and add the weak form for the incompressibility constraint we obtain the weak formulation:
\begin{equation}\begin{cases}\begin{aligned}
& \mbox{Find } (u_h,p_h) \in \mathcal{V}_h(\Omega)\times \mathcal{Q}_h(\Omega) \mbox{ such that:} \\
& \phantom{+} \int_\Omega \frac12 \left( v \otimes w : \nabla u - u \otimes w : \nabla v\right) {\rm d}V \\
& + \int_{\Gamma^D} \frac12 {\rm max}\left( 0, n \cdot g^D \right) \left( v \cdot u \right) {\rm d}S + \int_{\Gamma^N} \frac12 {\rm max}\left( 0, n \cdot w\right) \left( v \cdot u \right) {\rm d}S \\
& + 2 \nu \int_\Omega \nabla^s v : \nabla^s u {\rm d}V + 2 \nu \int_{\Gamma^D} - vn : \nabla^s u - un : \nabla^s v + \beta \left(v\cdot u\right) {\rm d}S \\
& + \int_\Omega  - p {\rm div}\left( v \right) {\rm d}V + \int_{\Gamma^D} p \left( v \cdot n \right) {\rm d}S \\
= & \phantom{+} \int_{\Gamma^D} - \frac12 {\rm min} \left( 0, n \cdot g^D \right) \left( v \cdot g^D \right) {\rm d}S \\
 & + 2 \nu \int_{\Gamma^D} - g^D n : \nabla^s v + \beta \left(v\cdot g^D\right) {\rm d}S \hspace{5.15cm} \forall v \in \mathcal{V}_h(\Omega), \\\\
& \phantom{+} \int_\Omega - q {\rm div}\left( u \right) {\rm d}V  + \int_{\Gamma^D} q \left( u \cdot n \right) {\rm d}S \\
= & \phantom{+} \int_{\Gamma^D} q \left( g^D \cdot n \right) {\rm d}S \hspace{8.05cm} \forall q \in \mathcal{Q}_h(\Omega).
\end{aligned}\end{cases}\label{eq:NSfull}\end{equation}
Similar formulations can be found in \cite{Burman2007,Burman2009,Bazilevs2012,Burman2014,Hsu2014,Schott2014,Kamensky2015,Schott2015,Xu2015,Hsu2016,Schott2016,Villanueva2017,Wu2017}. The variational forms in most of these references contain additional volumetric stabilization terms to enhance stability for large convective velocities, similar to the SUPG stabilization terms we added for the convection-diffusion problem in Section~\ref{sec:SUPG}. The numerical examples with ${\rm Re} \approx 100$ presented here and with ${\rm Re} \approx 2000$ presented in Section~\ref{sec:NSvalidation} did not require these, however. Furthermore, SUPG (or more general variational multiscale (VMS) \cite{Hughes2000,Hughes2004,Bazilevs2007Variational}) stabilization terms have originally been developed for linear bases and have in the context of immersed finite element methods only been used with such linear bases \cite{Bazilevs2012,Xu2015,Hsu2016}, with a stabilization parameter that is decreased in the vicinity of the trimmed elements \cite{Hsu2014,Kamensky2015,Wu2017} or in combination with additional stabilization techniques such as ghost penalty terms (see \emph{e.g.,}\ \cite{Burman2010,BurmanHansbo2012}) on trimmed elements \cite{Schott2014,Schott2015,Schott2016,Kadapa2017,Villanueva2017}.

\begin{figure}
 \centering
 \begin{subfigure}[t]{.49\textwidth}
  \includegraphics[width=\textwidth]{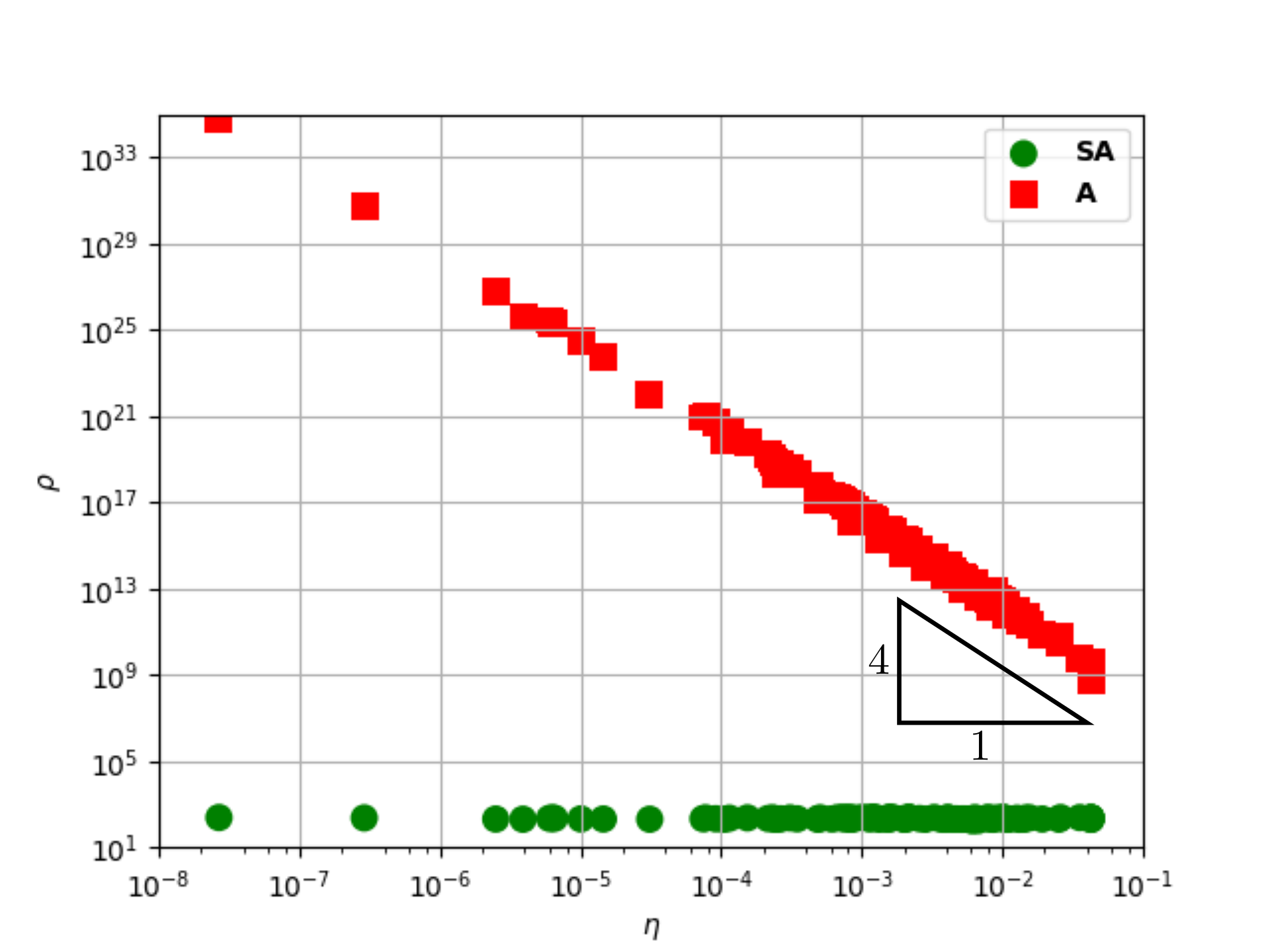}
  \caption{}
  \label{fig:NS_eta}
 \end{subfigure}
 \hfill
 \begin{subfigure}[t]{.49\textwidth}
  \includegraphics[width=\textwidth]{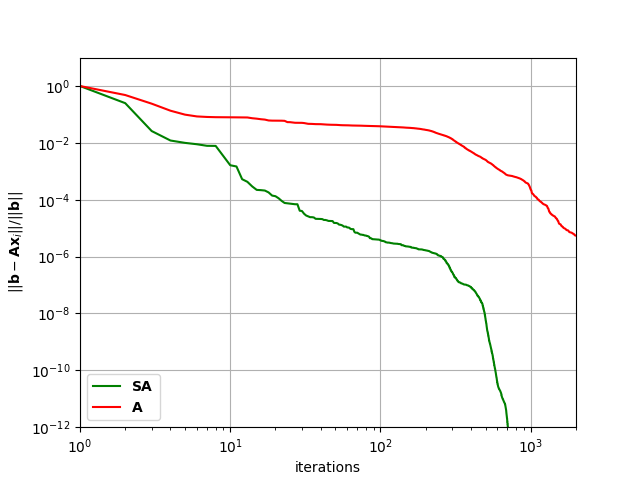}
  \caption{}
  \label{fig:NS_residual}
 \end{subfigure}
 \caption{(a) Original and preconditioned eigenvalue ratios for the Navier-Stokes problem. (b) The relative (preconditioned) residual error versus the number of GMRES iterations for an arrangement with $\theta = \SI{25}{\degree}$ yielding $\eta=9\cdot10^{-4}$. \label{fig:NS}}
\end{figure}

The system matrix and consequently the conditioning depend on the convective velocity $w$. As every Picard iteration uses a different convective velocity $w$, the conditioning is different for every Picard iteration. In order to keep the results clear, we only present the eigenvalue ratios with the converged convective velocity. We have observed that the convective velocity does not significantly affect the conditioning, as also follows from the results in Section~\ref{sec:NSvalidation}. Similar to the Stokes problem, the matrices are preconditioned by the CbAS preconditioner using the Schur complement. The eigenvalue ratios with and without preconditioning are presented in Figure~\ref{fig:NS_eta}. Comparison to Figure~\ref{fig:Stokes_eta} conveys that these eigenvalue ratios are essentially identical to those for the stokes case. Therefore the same observations apply, \emph{i.e.,}\ the eigenvalue ratios without preconditioning scale with the smallest volume fraction to the power $4$ and the CbAS preconditioner provides a well-conditioned system that is robust with respect to the smallest volume fraction. All obtained eigenvalue ratios are within the range of $170 \leq \rho(\mathbf{A}) \leq 244$. A typical GMRES converges plot that illustrates the improved convergence behavior using CbAS is shown in Figure~\ref{fig:NS_residual}.

\subsubsection{Transient Navier-Stokes problem}\label{sec:NSvalidation}
We also evaluate the conditioning of the Navier-Stokes system for a transient problem. We set the dynamic viscosity to $\nu = 10^{-3}$ and enlarge the domain according to:
\begin{equation}
 \Omega = \bigg\{ x = \big(x_1,x_2\big) \bigg| x_1 \in \big(-2,6\big) \hspace{.1cm} , \hspace{.1cm} x_2 \in \big(-1,1\big) \hspace{.1cm} , \hspace{.1cm} \big|x+c\big|_e > \frac14 \bigg\}.
\end{equation}
With a Poiseuille inflow profile with maximal velocity $1$ this yields a Reynold's number of approximately $2000$, which is expected to result in an unsteady flow. To perturb symmetry in the solution we consider a slightly off-centered circular inclusion, \emph{i.e.,}\ $c = (0,10^{-3})$. On the lower and upper boundaries ($|x_2|=1$) and on the boundary of the circular exclusion ($|x+c|_e=\frac14$) we prescribe homogeneous Dirichlet conditions, \emph{i.e.,}\ no-slip, and on the right boundary ($x_1=6$) we again prescribe the homogeneous directional do-nothing outflow condition. The Poiseuille inflow condition reads $u = (1-x_2^2,0)^T$. For the initial condition we use the solution of the Stokes problem posed on the same domain and subject to the same boundary conditions.

The size of the encapsulating domain is increased such that the same function space and grid size as used for the steady Navier-Stokes problem discussed above can be used. The angle between the grid lines and the domain is set to $\theta = \pi/4$. We use a time step $\delta t = 10^{-2}$ and employ the following $\vartheta$-scheme with $\vartheta = \frac12$ (Crank-Nicolson) for time integration:
\begin{equation}
 \underbrace{\begin{bmatrix}
  \mathbf{M} + \vartheta \delta t \mathbf{A}_{vu}^\vartheta & \delta t \mathbf{A}_{vp} \\ \delta t \mathbf{A}_{vp}^T & 
 \end{bmatrix}}_{=\mathbf{K}}
 \underbrace{\begin{bmatrix}
  \mathbf{x}_u^{n+1} \\
  \mathbf{x}_p^\vartheta
 \end{bmatrix}}_{=\mathbf{x}}
 =
 \underbrace{\begin{bmatrix}
  \left( \mathbf{M} - (1-\vartheta)\delta t \mathbf{A}_{vu}^\vartheta \right) \mathbf{x}_u^n + \delta t \mathbf{b}_v \\
  \delta t \mathbf{b}_q
 \end{bmatrix}}_{=\mathbf{b}}.
\end{equation}
The matrix $\mathbf{A}$ and right hand side vector $\mathbf{b}$ follow from \eqref{eq:NSfull}. $\mathbf{A}_{vu}^\vartheta$ depends on the convective velocity $u_h^\vartheta = \vartheta u_h^{n+1} + (1-\vartheta)u_h^n$. Matrix $\mathbf{M}$ denotes the mass matrix. To solve this nonlinear system we employ the same Picard iterations as for the steady Navier-Stokes problem at every time step. For the initial guess we use the velocity $u_h^n$ at the previous time step. 

\begin{figure}
 \centering
 \begin{subfigure}[t]{.49\textwidth}
  \includegraphics[width=\textwidth]{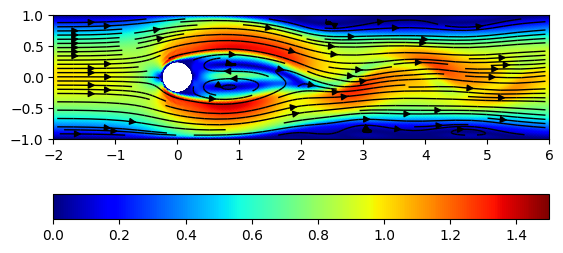}
  \caption{Velocity magnitude at $t=10$.}
 \end{subfigure}
 \hfill
 \begin{subfigure}[t]{.49\textwidth}
  \includegraphics[width=\textwidth]{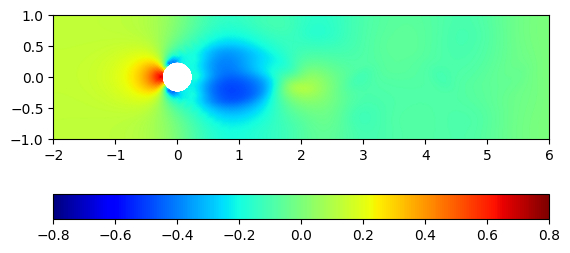}
  \caption{Pressure at $t=10$.}
 \end{subfigure}
 \\
 \begin{subfigure}[t]{.49\textwidth}
  \includegraphics[width=\textwidth]{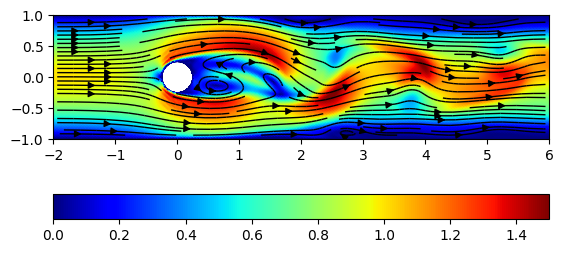}
  \caption{Velocity magnitude at $t=15$.}
 \end{subfigure}
 \hfill
 \begin{subfigure}[t]{.49\textwidth}
  \includegraphics[width=\textwidth]{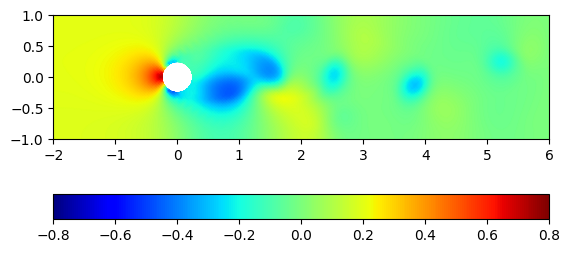}
  \caption{Pressure at $t=15$.}
 \end{subfigure}
 \\
 \begin{subfigure}[t]{.49\textwidth}
  \includegraphics[width=\textwidth]{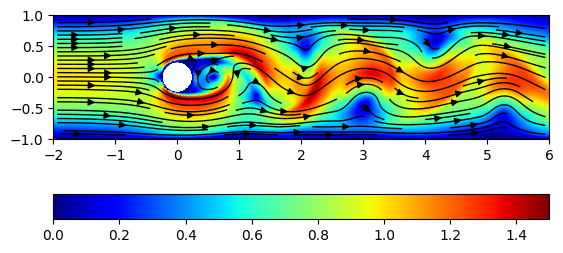}
  \caption{Velocity magnitude at $t=20$.}
 \end{subfigure}
 \hfill
 \begin{subfigure}[t]{.49\textwidth}
  \includegraphics[width=\textwidth]{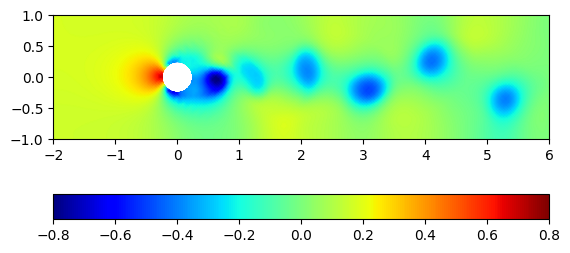}
  \caption{Pressure at $t=20$.}
 \end{subfigure}
 \\
 \begin{subfigure}[t]{.49\textwidth}
  \includegraphics[width=\textwidth]{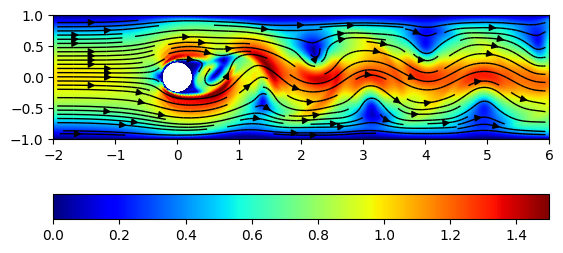}
  \caption{Velocity magnitude at $t=25$.}
 \end{subfigure}
 \hfill
 \begin{subfigure}[t]{.49\textwidth}
  \includegraphics[width=\textwidth]{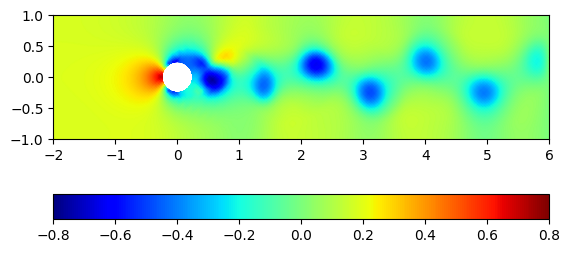}
  \caption{Pressure at $t=25$.}
 \end{subfigure}
 \caption{Various snapshots of flow profiles of the weakly turbulent immersed Navier-Stokes computation. \label{fig:NStransient}}
\end{figure}

The resulting flow profiles are presented in Figure~\ref{fig:NStransient}. The flow displays unsteady behavior and after a transient regime it settles into a periodic Von K\'{a}rm\'{a}n vortex shedding for $t \gtrsim 2$. The evolution of the original and preconditioned eigenvalue ratios of matrix $\mathbf{K}$ at the final Picard

\begin{figure}
 \centering
  \includegraphics[width=.7\textwidth]{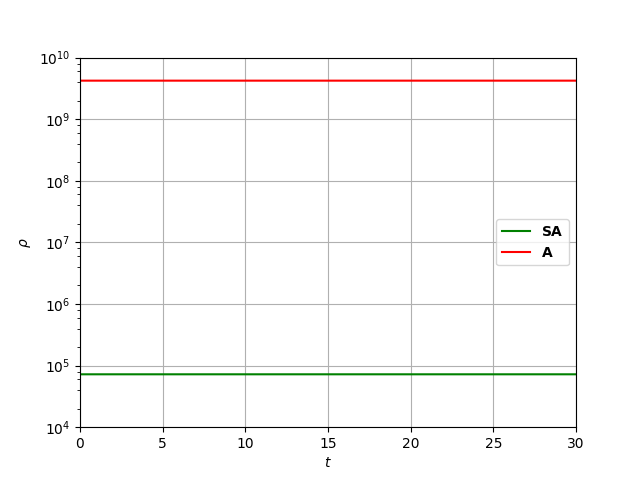}
  \caption{Original and preconditioned eigenvalue ratios for the weakly turbulent immersed Navier-Stokes computation plotted against time. The preconditioner significantly improves the conditioning of the linear system and both the original and preconditioned eigenvalue ratios are virtually independent of the convective velocity.\label{fig:NStransientConditioning}}
\end{figure}

\noindent iteration are presented in Figure~\ref{fig:NStransientConditioning}. From the fact that the eigenvalue ratios are essentially time-independent, we infer that the convective velocity at this Reynold's number does not significantly affect the conditioning. It is observed that the CbAS preconditioned system is well-conditioned.

\section{Concluding remarks}\label{sec:conclusion}
A challenging aspect of immersed methods is ill-conditioning, which commonly precludes the computation of high-fidelity solutions. The SIPIC preconditioning technique developed in \cite{SIPIC} effectively resolves this problem for symmetric positive definite (SPD) systems with smooth isogeometric bases. However, it is inapplicable to problems that do not fall within this restricted class. Most notably, SIPIC does not provide a solution for preconditioning immersed finite element computations of flow problems. The main goal of this work has been to generalize the SIPIC preconditioner to immersed problems that are not SPD, and in particular to enable the preconditioning of immersed flow problems, and to improve the robustness for non-smooth bases.

The root of the ill-conditioning problem for systems that are not SPD is not different from that of SPD systems, \emph{viz.}\ on cut elements with small volume fractions  (\emph{i.e.,}\ the fraction of an element that intersects the physical domain) basis functions can become excessively small or almost linearly dependent. Consequently, we observe scaling relations between the conditioning and the smallest volume fraction for systems that are not SPD similar to those for SPD systems.

The Connectivity-based Additive-Schwarz (CbAS) preconditioner developed herein enables consideration of non-SPD single-field systems. Degrees of freedom that are supported on trimmed elements are inverted block-wise in accordance with the well-established Additive-Schwarz framework. This assembly procedure only requires the corresponding matrix blocks to be invertible, thereby relieving the constraint that the system matrix must be SPD. The essential aspect of the CbAS preconditioner is the way in which the blocks are constructed. To satisfy the condition that all almost linearly dependent basis functions are present within a block, basis functions are grouped based on their element connectivity. This makes the technique robust for high-order bases with low regularity. The CbAS preconditioner is essentially equivalent to SIPIC for SPD systems with spline bases with full regularity, but generalizes it to non-SPD single-field systems.

Efficient preconditioning of mixed formulations requires a field-wise treatment of the CbAS preconditioner. For mixed forms such as encountered in the Stokes and Navier-Stokes equations, we presented a field-wise preconditioning procedure. In the case of a mixed system of the form of a \mbox{(Navier-)Stokes} problem, the single-field CbAS preconditioner is first applied to the velocity-velocity block of the system for each spatial direction, and subsequently the Schur complement is used to compute an optimal CbAS preconditioner for the velocity-pressure block.

We have performed an extensive set of simulations to test the performance of the CbAS preconditioner, which demonstrate the effectiveness of the preconditioner for nonsymmetric single-field and mixed problems, and for immersed flow problems. Most importantly, for all considered cases, the conditioning of the preconditioned matrix is independent of the smallest trimmed-element volume fraction. The presented results convey that the CbAS preconditioner enables the computation of high-fidelity solutions to immersed fluid-flow problems using iterative solvers.

In our simulations we have restricted ourselves to uniform meshes and moderate Reynold's numbers. Although the CbAS preconditioning technique yielded near optimal results for uniform meshes, this does not necessarily extend to non-uniform grids with \emph{e.g.,}\ local or hierarchical refinements since the connectivity structure can in such cases be fundamentally different from the cases considered here. Moreover, application of CbAS to higher Reynold's number flows, possibly with additional stabilization terms, could require more advanced domain decomposition or multigrid techniques such as coarse grid corrections.

\section*{Acknowledgement}
\noindent The research of F.\ de Prenter was funded by the NWO under the Graduate Program Fluid \& Solid Mechanics. All simulations in this work were performed using the open source software package Nutils (www.nutils.org). We would like to acknowledge Christoph Lehrenfeld for fruitful discussions regarding Additive-Schwarz preconditioning.

\section*{References}
\bibliographystyle{elsarticle-num}
\bibliography{RefFileFrits}

\end{document}